\documentclass[11pt,leqno]{article}

\usepackage{amsmath,amsfonts,amscd,amssymb,theorem}

\long\def\comment#1\endcomment{}

\comment
\endcomment

\makeatletter
\begingroup
\gdef\th@dotted{\normalfont\itshape
  \def\@begintheorem##1##2{%
        \item[\hskip\labelsep \theorem@headerfont ##1\ ##2.]}%
\def\@opargbegintheorem##1##2##3{%
   \item[\hskip\labelsep \theorem@headerfont ##1\ ##2\ (##3).]}}
\endgroup
\makeatother

\theoremstyle{dotted}

\newtheorem{theorem}{Theorem}[section]
\newtheorem{lemma}[theorem]{Lemma}

\newtheorem{prop}[theorem]{Proposition}
\newtheorem{corr}[theorem]{Corollary}

\makeatletter
\begingroup
\gdef\th@upshape{\normalfont
  \def\@begintheorem##1##2{%
        \item[\hskip\labelsep \theorem@headerfont ##1\ ##2.]}%
\def\@opargbegintheorem##1##2##3{%
   \item[\hskip\labelsep \theorem@headerfont ##1\ ##2\ (##3).]}}
\endgroup
\makeatother

\theoremstyle{upshape}

\newtheorem{defn}[theorem]{Definition}
\newtheorem{remark}[theorem]{Remark}

\makeatletter
\renewcommand{\subsection}{\@startsection{subsection}{2}{0pt}{-3ex
plus -1ex minus -0.2ex}{-2mm plus -0pt minus
-2pt}{\normalfont\bfseries}} 
\renewcommand{\subsubsection}{\@startsection{subsubsection}{3}{0pt}{-3ex
plus -1ex minus -0.2ex}{-2mm plus -0pt minus
-2pt}{\normalfont\bfseries}} 
\makeatother

\makeatletter
\@addtoreset{equation}{section}
\makeatother

\newcommand{\cntrct}                
{\hspace{2pt}\raisebox{1pt}{\text{$\lrcorner$}}\hspace{2pt}}

\newcommand{\proof}[1][Proof.]{\smallskip\noindent{\em #1}}
\def\endproof{\hfill\ensuremath{\square}\par\medskip}

\renewcommand{\labelenumi}{{\normalfont(\roman{enumi})}}

\def\eqref#1{\thetag{\ref{#1}}}

\let\latexref=\ref
\def\ref#1{{\normalfont{\latexref{#1}}}}

\newcommand{\wt}{\widetilde}

\newcommand{\ol}{\overline}

\newcommand{\dg}{\dagger}

\newcommand{\idot}{{\:\raisebox{1pt}{\text{\circle*{1.5}}}}}
\newcommand{\hdot}{{\:\raisebox{3pt}{\text{\circle*{1.5}}}}}
\newcommand{\eps}{\varepsilon}
\renewcommand{\phi}{\varphi}

\newcommand{\vH}{\check{H}}

\newcommand{\Hom}{\operatorname{Hom}}
\newcommand{\Ext}{\operatorname{Ext}}

\newcommand{\Coker}{\operatorname{Coker}}
\newcommand{\Ker}{\operatorname{Ker}}
\renewcommand{\Im}{\operatorname{Im}}

\newcommand{\id}{\operatorname{\sf id}}

\newcommand{\gr}{\operatorname{\sf gr}}
\newcommand{\tr}{\operatorname{\sf tr}}
\newcommand{\colim}{\operatorname{colim}}

\newcommand{\D}{{\cal D}}
\newcommand{\C}{{\cal C}}

\newcommand{\K}{{\mathbb K}}

\newcommand{\hush}{\natural}
\newcommand{\hash}{\sharp}

\newcommand{\Sets}{\operatorname{Sets}}
\newcommand{\Cat}{\operatorname{Cat}}
\newcommand{\Aut}{{\operatorname{Aut}}}

\newcommand{\amod}{{\text{\rm -mod}}}

\newcommand{\ppt}{{\sf pt}}

\newcommand{\lotimes}{\overset{\sf\scriptscriptstyle L}{\otimes}}

\newcommand{\Spec}{\operatorname{Spec}}

\newcommand{\cchar}{\operatorname{\sf char}}

\newcommand{\Z}{{\mathbb Z}}
\newcommand{\N}{{\mathbb N}}

\newcommand{\ux}{\underline{x}}

\newcommand{\LR}{\Lambda R}

\newcommand{\Fun}{\operatorname{Fun}}

\newcommand{\e}{\operatorname{\sf e}}

\newcommand{\E}{\mathcal{E}}

\newcommand{\vpi}{\check{\pi}}

\newcommand{\Stab}{\operatorname{Stab}}

\newcommand{\bi}{\overline{i}}

\newcommand{\bW}{\overline{W}}

\title{Hochschild-Witt complex}

\author{D. Kaledin\thanks{Partially supported by RScF, grant
    14-21-00053}}

\begin{document}

\maketitle

\tableofcontents

\section*{Introduction.}

This paper concludes the project started in \cite{witt}. The goal of
the project is to find a non-commutative generalization of the
classical construction of the $p$-typical Witt vectors ring $W(A)$
associated to a commutative associative ring $A$. As we have
explained in the introduction to \cite{witt}, such a generalization
should also produce a generalization of the de Rham-Witt complex of
Deligne and Illusie \cite{ill}, and a generalization of the Witt
vectors group $W(A)$ constructed for any associative ring $A$ by
Hesselholt \cite{hewi}, \cite{heerr}.

The common theme unifying \cite{ill} and \cite{hewi} is Hochschild
homology. The relation between de Rham-Witt complex and Hochschild
homology is implicit, and it goes through the identification of
differential forms and Hochschild homology classes discovered by
Hochschild, Kostant and Rosenberg \cite{HKR}. To take account of the
de Rham differential, one has to use cyclic homology (\cite{Lo},
\cite{FTadd}). The relation between Hesselholt's construction and
Hochschild homology is much more direct, and it goes back to the
theory of Topological Cyclic Homology and cyclotomic trace of
B\"okstedt, Hsiang and Madsen (\cite{BHM}, see also an exposition in
\cite{HM}). As a part of their theory, B\"okstedt, Hsiand and Madsen
construct a spectrum $TR(A;p)$ for any ring spectrum $A$, and in
particular, for any associative ring $A$. This spectrum is
manifestly related to Hochschild homology --- in fact, it is
constructed starting from B\"okstedt's Topological Hochschild
Homology spectrum $THH(A)$ of \cite{bo}. What Hesselholt did in
\cite{hewi} was the following: he took an arbitrary associative ring
$A$, and gave a purely algebraic construction of the homotopy group
$\pi_0TR(A;p)$.

It is known that if the ring $A$ is in fact an algebra over a field
$k$ of characteristic $p$ --- which happens to be the most
interesting case for the Witt vector construction --- then $TR(A;p)$
is an Eilenberg-Mac Lane spectrum. Thus effectively, $TR(A;p)$ is a
chain complex $TR_\idot(A;p)$, and its homotopy groups $\pi_\idot
TR(A;p)$ are the homology groups of the complex $TR_\idot(A;p)$. In
the situation of the Hochschild-Kostant-Rosenberg Theorem --- that
is, if $A$ is commutative and smooth --- it has been further proved
by Hesselholt \cite{hedRW} that these homology groups are naturally
identified with the terms of the de Rham-Witt complex of $X=\Spec
A$. This established the link between \cite{ill} and \cite{BHM} and
became the subject of a lot of further research.

For a general associative unital $k$-algebra $A$, one still has the
complex $TR_\idot(A;p)$, but there is no ready algebraically defined
complex one could compare it with. A natural goal is, then, to try
to construct such a complex.

\medskip

This is not what we do in this paper, or at least, not quite. What
we do is construct a functorial ``Hochschild-Witt complex''
$WCH_\idot(A)$ for any associative unital algebra over a perfect
field $k$ of characteristic $p$ that looks very much like
$TR_\idot(A;p)$ and carries all the additional structures that
$TR(A;p)$ is known to carry. However, our construction is purely
homological, and at present, we do not know how to compare it to
equivarient stable homotopy theory used in \cite{BHM}. Lacking a
comparison theorem between $WHH_\idot(A)$ and $TR_\idot(A;p)$, we at
least prove that in degree $0$, the homology group $WHH_0(A)$ of the
complex $WCH_\idot(A)$ is indeed Hesselholt's Witt vectors group,
and in the Hochschild-Kostant-Rosenberg situation, our construction
recovers the de Rham-Witt complex of Deligne and Illusie.

\medskip

As explained in the introduction to \cite{witt}, our main method is
to observe that Hochschild homology is in fact a theory with two
variables, an algebra $A$ and an $A$-bimodule $M$. Thus it is
natural to expect to have a two-variable ``Hochschild-Witt
homology'' functors $WHH_\idot(A,M)$. Moreover, if one can manage to
equip with functors with an additional structure of a ``trace
theory'' in the sense of \cite{trace}, then one only needs to
construct $WHH_\idot(A,M)$ for $A=k$, the base field. This is
exactly what has been done in \cite{witt}. After that, the groups
$WHH_\idot(A,M)$ can be simply produced by the general machine of
\cite{trace}, and taking $M=A$, the diagonal bimodule, one recovers
the homology groups $WHH_\idot(A) = WHH_\idot(A,A)$ of the
Hochschild-Witt complex $WCH_\idot(A)$.

However, the general machine gives very little information on the
groups $WHH_\idot(A)$, nor allows one to prove any comparison
theorems. From this point of view, the Hochschild-Witt complex
$WCH_\idot(A)$ deserves a separate detailed study. This is exactly
what the present paper aims to provide.

\medskip

Let us give a brief overview of the paper. We start with
preliminaries, and there is a lot of them --- to the point of taking
up whole two sections, Section~\ref{prem1.sec} and
Section~\ref{prem2.sec}. The reason for this is that our approach to
Hochschild and cyclic homology is based on the fundamental principle
discovered by Connes \cite{connes}: a large part of the formalism of
cyclic homology has nothing to do with algebras and should be
developed for arbitrary functors from the cyclic category $\Lambda$
to vector spaces or abelian groups. Both Section~\ref{prem1.sec} and
Section~\ref{prem2.sec} work in this generality, thus deal
essentially with linear algebra.

The material in Section~\ref{prem1.sec} is completely standard. In
fact, the only statements for which we do not have ready references
are Corollary~\ref{i.l.triv} and Lemma~\ref{fp.le}, and even these
are too simple to be new. Everything else could be replaced with a
reference to any of the standard sources such as \cite{Lo} or
\cite{FTadd}. However, the results are there not quite in the form
we need, the assembly is different, and different notation would
require a lot of translation. In the end, we have decided that it
would be cheaper and faster to develop the theory from scratch. The
side effect of this is that as far as cyclic homology is concerned,
the paper is more-or-less self-contained.

Section~\ref{prem2.sec} still deals with linear algebra but the
facts we need are less standard; they are mostly concerned with
cyclic covers $\Lambda_l$ of the cyclic category $\Lambda$ that go
back to \cite[Appendix]{FTadd}. We do not have ready references for
these results, and some of them (e.g.\
Proposition~\ref{Y.prop}) do look new.

Section~\ref{cons.sec} is still essentially concerned with
preliminaries. In Subsection~\ref{alg.subs}, we recall the
construction of the cyclic homology of an associative algebra and
some basic comparison theorems. In Subsection~\ref{car.subs}, we
recall a non-commutative generalization of the Cartier isomorphism
developed in \cite{cart}. Subsection~\ref{trace.subs} recapitulates
the machinery of trace functors of \cite{trace}, and
Subsection~\ref{witt.subs} recalls the basic properties of the
polynomial Witt vectors functors constructed in \cite{witt}.

After all these preliminaries, Section~\ref{defs.sec} finally
introduces the main object of our study, the Hochschild-Witt complex
$WCH_\idot(A)$ and the Hochschild-Witt homology groups
$WHH_\idot(A)$ of an associative unital algebra $A$ over a perfect
field $k$ of positive characteristic $p$. We do this in
Subsection~\ref{basic.subs}, and we also construct some additional
structures on these groups --- namely, the Verschiebung and
Frobenius maps $V$, $F$, and the Teichm\"uller representative map
$T$. We prove all sorts of identities these maps satisfy. Then in
Subsection~\ref{free.subs}, we illustrate the general theory by
computing explicitly the Hochschild-Witt homology groups of a tensor
algebra $T^\hdot(M)$ of a $k$-vector space $M$.

The last two sections contain our comparison
theorems. Section~\ref{hess.sec} proves that in degree $0$, we
recover Witt vectors group constructed by Hesselholt. In
Section~\ref{ill.sec}, we prove that in the
Hochschild-Kostant-Rosenberg situation, our Hochschild-Witt homology
groups are naturally identified with the terms of the de Rham-Witt
complex of \cite{ill}.

Let us mention that even in the commutative smooth case, where we
effectively just give a new construction of the de Rham-Witt
complex, we believe that our approach is instructive. It is
significantly less {\em ad hoc} than \cite{ill}, and it explains
some mysteries of the theory. One such is the structure of the
associated graded quotient of the de Rham-Witt complex with respect
to the natural filtration. While in principle, it was completely
described by Illusie, we believe that our general non-commutative
expression given in Corollary~\ref{gr.corr} is much simpler, and it
is not possible to really understand what is going on without using
the cyclic homology approach.

Two comparison results are conspicuously missing from the present
paper. One is comparison with $TR(A;p)$ mentioned earlier; as we
said, we expect it to be true, but we currently do not have enough
technology to prove it. Another is a very natural lifting result
that is very important in \cite{ill} --- namely, if our $k$-algebra
$A$ lifts to a flat algebra $\wt{A}$ over the Witt vectors ring
$W(k)$, then the periodic cyclic version of our de Rham-Witt
homology should coincide with the periodic cyclic homology of
$\wt{A}$. Unfortunately, we could not prove it, and at present it is
not even clear if it is true. If it is true, there is further
ambiguity as to what the precise statement should be --- instead of
the usual periodic cyclic homology, one might need to use the
co-periodic cyclic homology introduced recently in \cite{coper}. All
in all, this certainly deserves further research.

Yet another very intriguing object for comparison is the recent very
general construction of Witt vectors and de Rham-Witt complex given
by Cuntz and Deninger \cite{CD}. Among other things, they managed to
define a Witt vectors group for an arbitrary associative ring $A$ in
such a way that it is also an associative ring. This goes against
the main trend of our paper since Hochschild homology $HH_0(A)$ is
only a ring for a commutative $A$, so that the relation between our
construction and \cite{CD} cannot be straightforward. At this point,
we do not know what the relation might be.

\subsection*{Acknowledgements.}

This paper as well as \cite{witt} took a long time to appear; in
fact, the results were announced more than five years ago. I
apologize for the delay, and I am very grateful to many people with
whom I had the opportunity to discuss the subject. As with
\cite{witt}, I am especially grateful to L. Hesselholt whose
generous input was absolutely crucial for most of this work, and to
V. Vologodsky, to whom the main technical idea of \cite{witt} is
largely due.

\section{Preliminaries I.}\label{prem1.sec}

\subsection{Generalities.}

We follow the same notations and conventions as \cite{witt}. In
particular, for any category $\C$, we denote by $\C^o$ the opposite
category, and for any two objects $c,c' \in \C$, we denote by
$\C(c,c')$ the set of morphisms from $c$ to $c'$. We denote by
$\ppt$ the point category, and for any group $G$, we denote by
$\ppt_G$ the groupoid with one object with automorphism group
$G$. For any integer $l \geq 1$, we simplify notation by letting
$\ppt_l = \ppt_{\Z/l\Z}$.

For any small category $I$ and any category $\C$, we will denote by
$\Fun(I,\C)$ the category of functors from $I$ to $\C$. For any ring
$A$, we denote by $A\amod$ the category of left $A$-modules, and for
any small category $I$, we will simplify notation by letting
$\Fun(I,A) = \Fun(I,A\amod)$. This is an abelian category; we denote
by $\D(I,A)$ its derived category. For any functor $\gamma:I \to I'$
between small categories, we denote by $\gamma^*:\Fun(I',A) \to
\Fun(I,A)$ the pullback functor, and we denote by
$\gamma_!,\gamma_*:\Fun(I,A) \to \Fun(I',A)$ its left and right
adjoint (the left and right Kan extensions along $\gamma$). The
functor $\gamma^*$ is exact, hence descends to derived categories,
and the derived functors $L^\hdot\gamma_!,R^\hdot\gamma_*:\D(I,A)
\to \D(I',A)$ are left and right adjoint to $\gamma^*:\D(I',A) \to
\D(I,A)$. The homology and the cohomology of a small category $I$
with coefficients in a functor $E \in \Fun(I,A)$ are given by
$$
H_\idot(I,E) = L^\hdot\tau_!E, \qquad H^\hdot(I,E) = R^\hdot\tau_*E,
$$
where $\tau:I \to \ppt$ is the tautological projection. If the ring
$A$ is commutative, then the category $\Fun(I,A)$ is a symmetric
unital tensor category with respect to the pointwise tensor
product. Its cohomology $H^\hdot(I,A)$ with coefficients in the
constant functor $A$ is then an algebra, and for any $E \in
\Fun(I,A)$, both $H_\idot(I,E)$ and $H^\hdot(I,E)$ are modules over
$H^\hdot(I,A)$.

We recall that if $I = \Delta^o$ is the opposite of the category
$\Delta$ of finite non-empty totally ordered sets, then functors $E
\in \Fun(\Delta^o,A)$ are simplicial $A$-modules, and homology
groups $HH_\idot(\Delta^o,E)$ can be computed by the {\em standard
  chain complex} $C_\idot(E)$ of the simplicial $A$-module $A$. The
term $C_i(E)$ of the complex $C_\idot(E)$ is the value of $E$ at the
set $\{0,\dots,i\}$ with the standard total order, and the
differential is the alternating sum of the face maps.

We will assume known the notions of a fibration, a cofibration and a
bifibration of small categories originally introduced in \cite{sga},
together with the correspondence between fibrations and
pseudofunctors known as the Grothendieck construction. A fibration,
cofibration or bifibration is {\em discrete} if its fibers are
discrete categories (that is, the only morphisms are identity
maps). We also assume known the following useful base change lemma:
if we are given a cartesian square
$$
\begin{CD}
I'_1 @>{f_1}>> I_1\\
@V{\pi'}VV @VV{\pi}V\\
I' @>{f}>> I
\end{CD}
$$
of small categories, and $\pi$ is a cofibration, then $\pi'$ is a
cofibration, and the base change map $L^\hdot\pi'_! \circ f_1^* \to
f^* \circ L^\hdot\pi_!$ is an isomorphism. Dually, if $\pi$ is a
fibration, then $\pi'$ is a fibration, and $f^* \circ R^\hdot\pi_*
\cong R^\hdot\pi'_* \circ f_1^*$. For a proof, see e.g.\ \cite[Lemma
  1.7]{ka0}. Another result contained in \cite[Lemma 1.7]{ka0} is the
following projection formula: for any cofibration $\gamma:I' \to I$
and $E \in \Fun(I,A)$, we have a natural isomorphism
\begin{equation}\label{proj.eq}
L^\hdot\pi_!\pi^*E \cong E \otimes_\Z L^\hdot\pi_!\Z.
\end{equation}
One specific example of a base change situation is a bifibration
$\pi:I' \to I$ whose fiber is equivalent to a groupoid $\ppt_G$. In
this case, for any $E \in \Fun(I',A)$ and any object $i' \in I'$
with image $i=\pi(i') \in I$, the $A$-module $E(i')$ carries a
natural action of the group $G$, and we have base change
isomorphisms
$$
\pi_!E(i) \cong E(i')_G, \qquad \pi_*E(i) \cong E(i')^G,
$$
where in the right-hand side, we have coinvariants and invariants
with respect to $G$. If the group $G$ is finite, then for any
$A[G]$-module $V$, we have a natural trace map
\begin{equation}\label{tr.G}
\tr_G = \sum_{g \in G}g:V_G \to V^G
\end{equation}
whose cokernel is the Tate cohomology group $\vH^0(G,V)$. Taken
together, these maps then define a natural trace map
\begin{equation}\label{tr.gamma}
\tr_\pi:\pi_!E \to \pi_*E.
\end{equation}
This map is functorial in $E$ and compatible with the base
change. In some cases, it is an isomorphism --- for example, this
happens if for any object $i' \in I'$, the $A[G]$-module $E(i')$ is
of the form $M[G]$ for some $A$-module $G$ (in such a situation, we
will say that $E$ is {\em $\pi$-free}). In the general case, we
denote by
\begin{equation}\label{v.ga}
\vpi_*:\Fun(I',A) \to \Fun(I,A)
\end{equation}
the functor sending $E$ to the cokernel of the map \eqref{tr.gamma}.
For any $i' \in I'$ with $i=\pi(i')$, we have a natural
identification $\vpi_*(E)(i) \cong \vH^0(G,E(i'))$.

We also note that for any $E \in \Fun(I',A)$, we can consider the
composition
$$
\begin{CD}
E @>{l}>> \pi^*\pi_!E @>{\pi^*(\tr_\pi)}>> \pi^*\pi_*E @>{r}>> E,
\end{CD}
$$
where $l$ and $r$ are the adjunction maps, and this defines a
functorial map
\begin{equation}\label{tr.dg}
\tr^\dg_\pi:E \to E.
\end{equation}
Moreover, the pullback functor $\pi^*:\Fun(I,A) \to \Fun(I',A)$ is
fully faithful. Therefore there exists a unique functorial map
\begin{equation}\label{e.gamma}
\e_\pi:\pi_*E \to \pi_!E
\end{equation}
such $\pi^*(\e_\pi)$ coincides with the composition
$$
\begin{CD}
\pi^*\pi_*E @>{r}>> E @>{l}>> \pi^*\pi_!E.
\end{CD}
$$

\subsection{Combinatorics.}\label{comb.subs}

There are many equivalent definitions of A. Connes' cyclic category
$\Lambda$ originally introduced in \cite{connes}; we will use the
one that embeds it into the category $\Cat$ of small categories. For
alternative approaches, see e.g.\ \cite[Chapter 6]{Lo},
\cite[Appendix]{FTadd}.

Denote by $[1]_\Lambda$ the category with one object whose
endomorphisms are non-negative integers $a \in \Z_+ \subset \Z$,
with composition given by sum. The embedding $\Z_+ \subset \Z$
induces a functor $[1]_\Lambda \to \ppt_{\Z}$. By the Grothendieck
construction, sets with a $\Z$-action correspond to small categories
discretely bifibered over $\ppt_\Z$, so that by pullback, every
$\Z$-set $S$ induces a discrete bifibration $[1]^S_\Lambda \to
[1]_\Lambda$. For any integer $n \geq 1$, we denote $[n]_\Lambda =
[1]^{\Z/n\Z}_\Lambda$ the category corresponding to the $\Z$-set
$\Z/n\Z$. Explicitly, objects in $[n]_\Lambda$ correspond to
residues $y \in \Z/n\Z$, and morphisms from $y$ to $y'$ are
non-negative integers $l$ such that $y+l=y' \mod n$. We denote the
set of objects of $[n]_\Lambda$ by $Y([n])$. Equivalently,
$[n]_\Lambda$ is the path category of the wheel quiver with $n$
vertices and $n$ edges, and $Y([n])$ is the set of vertices.

Any object $y \in Y([n])$ of the category $[n]_\Lambda$ has a
natural endomorphism $\tau_y:y \to y$ given by $l = n$, and for any
functor $f:[n]_\Lambda \to [m]_\Lambda$, we have
$$
f(\tau_y) = \tau_{f(y)}^{\deg f},
$$
where $\deg f$ is a non-negative integer independent of $y$, called
the {\em degree} of the functor $f$. The degree is obviously
multiplicative, $\deg(f_1 \circ f_2) = \deg(f_1)\deg(f_2)$.

\begin{defn}\label{v.h.def}
A functor $f:[n]_\Lambda \to [m]_\Lambda$ is {\em vertical} if it is
a discrete bifibration, {\em horizontal} if $\deg(f)=1$, and {\em
  non-degenerate} if $\deg(f) \neq 0$.
\end{defn}

\begin{defn}\label{LR.def}
The {\em cyclotomic category} $\LR$ is the small category whose
objects $[n]$ are numbered by positive integers $n \geq 1$, and
whose maps from $[n]$ to $[m]$ are given by non-degenerate functors
$f:[n]_\Lambda \to [m]_\Lambda$.
\end{defn}

Classes of horizontal resp.\ vertical maps are obviously closed
under compositions and contain all isomorphisms. Denote by
$\LR_h,\LR_v \subset \LR$ the subcategories with the same objects as
$\LR$ and horizontal resp.\ vertical maps between them. Then $\LR_v$
is naturally identified with the category of finite $\Z$-orbits ---
that is, finite sets equipped with a transitive $\Z$-action. On the
other hand, $\LR_h$ is the category $\Lambda$.

\begin{defn}
The {\em cyclic category $\Lambda$} is the subcategory $\Lambda =
\LR_h \subset \LR$.
\end{defn}

For any $[n] \in \LR$, the automorphism group $\Aut([n])$ of the
object $[n]$ is the cyclic group $\Z/n\Z$ (and all automorphisms are
both horizontal and vertical). Every functor $f:[n]_\Lambda \to
[m]_\Lambda$ of degree $1$ has left and right-adjoint functors
$f_\hash,f^\hash:[m]_\Lambda \to [n]_\Lambda$. Explicitly, they are
given by
\begin{equation}\label{lambda.du}
\begin{aligned}
f_\hash(y) &= \min\{ y' \in \Z | f(y') \geq \wt{y} \} \mod m, \\
f^\hash(y) &= \max\{ y' \in \Z | f(y') \leq \wt{y} \} \mod m,
\end{aligned}
\end{equation}
where we choose a representative $\wt{y} \in \Z$ of the residue
class $y$, and observe that the result does not depend on this
choice. Thus in particular, the category $\Lambda$ is equivalent to
its opposite category $\Lambda^o$, $\Lambda \cong \Lambda^o$.

\medskip

Denote by $\Lambda_\idot$ the category of vertical maps $v:[n] \to
[m]$ in $\LR$, with maps from $v:[n] \to [m]$ to $v':[n'] \to [m']$
given by commutative squares
$$
\begin{CD}
[n] @>{f}>> [n']\\
@V{v}VV @VV{v'}V\\
[m] @>{g}>> [m']
\end{CD}
$$
with horizontal $f$, $g$. Then $\Lambda_\idot$ splits as a disjoint
union
\begin{equation}\label{la.spl}
\Lambda_\idot = \coprod_{l \geq 1}\Lambda_l
\end{equation}
according to the degree $l = \deg(v)$ of the map $v$. The component
$\Lambda_1$ is the cyclic category $\Lambda$. Moreover, sending
$v:[n'] \to [n]$ to $[n']$ resp.\ $[n]$ gives two functors
\begin{equation}\label{i.pi.l}
i_\idot,\pi_\idot:\Lambda_\idot \to \Lambda
\end{equation}
and their individual components $i_l$, $\pi_l$, $l \geq 1$. The
functor $\pi_l$ is a bifibration with fiber $\ppt_l$. To simplify
notation, we denote an object $v:[n'] \to [n]$ of $\Lambda_l$ simply
by $[n]$; then automatically $n'=nl$, and we have $i_l([n]) = [nl]$,
$\pi_l([n]) = [n]$. For any object $[n] \in \Lambda_l$ represented
by a vertical arrow $v:[nl] \to [n]$, the arrow $v$ gives a natural
map
\begin{equation}\label{q.l}
\eta_l:Y(i_l([n])) \to Y(\pi_l([n])).
\end{equation}
The duality $\Lambda \cong \Lambda^o$ extends to
dualities $\Lambda_l \cong \Lambda_l^o$, $l \geq 1$.

For any two integers $m,l \geq 1$, a vertical map $v:[nml] \to [n]$
of degree $ml$ factorizes as
$$
\begin{CD}
[nml] @>{v'}>> [nm] @>{v''}>> [n],
\end{CD}
$$
with $v'$, $v''$ vertical of degrees $l$, $m$, and such a
factorization is unique up to a unique isomorphism. Sending $v$ to
the pair $\langle v',v'' \rangle$ then gives a commutative square
\begin{equation}\label{la.ml}
\begin{CD}
\Lambda_{ml} @>{i_m}>> \Lambda_l\\
@V{\pi_l}VV @VV{\pi_l}V\\
\Lambda_m @>{i_m}>> \Lambda,
\end{CD}
\end{equation}
and this square is obviously Cartesian. By abuse of notation, we
also denote by $\pi_l$, $i_m$ the left vertical and the top
horizontal functor in \eqref{la.ml} whenever there is no danger of
confusion.

For any $n,m,l \geq 1$, the functor $i_l$ induces a natural map
\begin{equation}\label{i.l.map}
\Lambda_l([n],[m]) \to \Lambda(i_l([n]),i_l([m])) =
\Lambda([ln],[lm]),
\end{equation}
and this map is obviously injective. We will need the following
stronger statement.

\begin{lemma}\label{fp.le}
Consider two objects $[n],[m] \in \Lambda$, and let $\Z/n\Z \times
\Z/m\Z = \Aut([n]) \times \Aut([m])$ act on $\Lambda([n],[m])$ by
$$
(g_1 \times g_2) \cdot f = g_2 \circ f \circ g_1^{-1}.
$$
Then the stabilizer $\Stab(f) \subset \Aut([n]) \times \Aut([m])$ of
any morphism $f \in \Lambda([n],[m])$ is a cyclic group $\Z/l\Z$
whose order divides both $n$ and $m$, its embedding into $\Aut([n])
\times \Aut([m])$ is the product of the standard embeddings $\Z/l\Z
\subset \Z/n\Z$, $\Z/l\Z \subset \Z/m\Z$, and the subset of fixed
points
$$
\Lambda([n],[m])^{\Stab(f)} \subset \Lambda([n],[m])
$$
coincides with the set $\Lambda_l([n/l],[m/l])$ embedded by the map
\eqref{i.l.map}.
\end{lemma}

\proof{} It is obvious from the definition that both individual
actions of $\Aut([n])$ and $\Aut([m])$ on $\Lambda([n],[m])$ are
free. Therefore when we project $\Stab(f) \subset \Aut([n]) \times
\Aut([m])$ to $\Aut([n])$ resp.\ $\Aut([m])$, the resulting maps
$\Stab(f) \to \Aut([n])$, $\Stab(f) \to \Aut([m])$ are
injective. This yields the first claim. Moreover, $f$ fits into a
commutative diagram
\begin{equation}\label{f.f.bar}
\begin{CD}
[n]_\Lambda @>{f}>> [m]_\Lambda\\
@VVV @VVV\\
[n]_\Lambda/\Stab(f) @>{\overline{f}}>> [m]_\Lambda/\Stab(f),
\end{CD}
\end{equation}
and if we denote by $\sigma \in \Stab(f)$ a generator of the cyclic
group $\Stab(f)\cong\Z/l\Z$, then we have $\sigma \circ f = f \circ
\sigma^{\deg{\overline{f}}}$. Since both vertical maps in
\eqref{f.f.bar} have degree $l$, we have
$\deg(\overline{f})=\deg(f)=1$, and this gives the second claim. The
third claim immediately follows from the definition of the functor
$i_l$.  \endproof

For every $[n] \in \Lambda$, a horizontal functor $h:[1]_\Lambda \to
[n]_\Lambda$ is uniquely determined by the image of the unique
object $o \in [1]_\Lambda$, so that we have a natural identification
$Y([n]) = \Lambda([1],[n])$. Moreover, for any $y \in Y([n])$,
composing the corresponding functor $h:[1]_\Lambda \to [n]_\Lambda$
with the unique functor $\ppt \to [1]_\Lambda$ gives gives an
embedding $i_y:\ppt \to [n]_\Lambda$ of the point category $\ppt$
onto the object $y$. For every map $f:[m] \to [n]$, we can define
the preimage $[f^{-1}(y)]$ by the fibered product
$$
\begin{CD}
[f^{-1}(y)] @>>> [n]_\Lambda\\
@VVV @VV{f}V\\
\ppt @>{i_y}>> [1]_\Lambda.
\end{CD}
$$
Then $[f^{-1}(y)]$ is the category corresponding to a finite totally
ordered set -- its objects are elements $y' \in f^{-1}(y) \subset
Y([m])$, and the category structure on $[f^{-1}(y)]$ induced from
$[n]_\Lambda$ defines a canonical total order on the set
$f^{-1}(y)$. In particular, take $n=1$. Then the category $\Delta$
of finite non-empty totally ordered sets is a full subcategory in
$\Cat$, and sending $f:[n] \to [1]$ to $[f^{-1}(o)]$ gives a functor
$$
\Lambda/[1] \to \Delta.
$$
This functor is an equivalence of categories. Choosing an inverse
equivalence and composing it with the forgetful functor, we obtain a
natural embedding
$$
j:\Delta \to \Lambda.
$$
More generally, for any horizontal map $h:[n] \to [1]$ and any
vertical map $v:[m] \to [1]$ in $\LR$, we have
$$
[n]_\Lambda \times_{[1]_\Lambda} [m]_\Lambda \cong [nm]_\Lambda,
$$
and sending $\langle h,v \rangle$ to their fibered product gives a
functor
\begin{equation}\label{wt.j}
\wt{j}:\Delta \times \LR_v \to \LR.
\end{equation}
Restricting $\wt{j}$ to $\Delta \times [l]$ for some $l \geq 1$, we
obtain natural embeddings
\begin{equation}\label{j.l}
\wt{j}_l:\Delta \times \ppt_l \to \LR, \qquad j_l:\Delta \to \Lambda_l,
\end{equation}
where we identify $\Z/l\Z \cong \Aut([l])$. By duality, \eqref{wt.j}
and \eqref{j.l} induce embeddings
$$
j^o:\Delta^o \to \Lambda, \qquad j^o_l:\Delta^o \to \Lambda_l, \quad
l \geq 1.
$$
The embeddings $j_l$ are compatible with the embeddings $i_l$ of
\eqref{i.pi.l} in the following sense: for any $l \geq 1$, we have a
commutative square
\begin{equation}\label{i.l.sq}
\begin{CD}
\Delta @>{j_l}>> \Lambda_l\\
@V{\bi_l}VV @VV{i_l}V\\
\Delta @>{j}>> \Lambda,
\end{CD}
\end{equation}
where the functor $\bi_l:\Delta \to \Delta$ is the functor sending a
totally ordered set $[n]$ with $n$ elements to the product $[l]
\times [n] \cong [nl]$ with lexicographical order. We have a similar
commutative square for the categories $\Delta^o$ and the functors
$j_l^o$. We note that if we denote by $s:[1] \to [l]$ the map that
sends the unique element in $[1]$ to the initial element in $[l]$,
then we have a functorial map
\begin{equation}\label{s.l}
s_l = s^o \times \id:\bi_l^o \to \id.
\end{equation}
The choice of $s$ breaks cyclic symmetry, so that an analogous map
for the functor $i_l$ does not exist.

\subsection{Homology.}\label{hom.subs}

Next, we recall briefly the main homological properties of the cyclic
categories $\Lambda_l$.

For any ring $A$ and $E \in \Fun(\Lambda,A)$, the functors
$j^o:\Delta^o \to \Lambda$, $j:\Delta \to \Lambda$ give canonical
objects $j^o_!j^{o*}E,j_*j^*E \in \Fun(\Lambda,A)$. Since $j$ is
equivalent to a discrete fibration and $j^o$ is equivalent to a
discrete cofibration, the Kan extensions functors $j_!^o$ and $j_*$
are exact by base change, and the projection formula \eqref{proj.eq}
gives canonical identifications
\begin{equation}\label{E.proj}
j^o_!j^{o*}E \cong E \otimes j^o_!\Z, \qquad j_*j^*E \cong E \otimes
j_*\Z,
\end{equation}
where $\Z$ is the constant functor with values $\Z$. Explicitly, one
has
\begin{equation}\label{k.exp}
j_!^o\Z([n]) \cong \Z[\Lambda([1],[n]), \qquad j_*\Z
  \cong \Z[\Lambda([n],[1])]
\end{equation}
for any $[n] \in \Lambda$, so that any map $f:[1] \to [n]$ gives a
element $[f] \in j^o_!\Z([n])$, and any map $f:[n] \to [1]$ gives an
element $[f] \in j_*\Z([n])$. Define a map $b:j_*\Z([n]) \to
j^o_!\Z([n])$ by
\begin{equation}\label{b.eq}
b([f]) = [f_\hash] - [f^\hash]
\end{equation}
for any $f:[n] \to [1]$, where $f^\hash$, $f_\hash$ are as in
\eqref{lambda.du}.

\begin{lemma}\label{4-le}
The maps $b$ of \eqref{b.eq} are functorial in $[n]$, and we have an
exact sequence
\begin{equation}\label{4.term}
\begin{CD}
0 @>>> \Z @>>> j_*\Z @>{b}>> j_!^o\Z @>>> \Z @>>> 0
\end{CD}
\end{equation}
of objects in $\Fun(\Lambda,\Z)$, where $\Z$ stands for the constant
functor.
\end{lemma}

\proof{} Identify $\Lambda([1],[n]) = Y([n])$ with the set $\Z/n\Z$
of objects in $[n]_\Lambda$, and identify $\Lambda([1],[n]) \cong
\Lambda([n],[1])$ by sending $f:[n] \to [1]$ to $f_\hash$. Then $b$
becomes a map $\Z[\Z/n\Z] \to \Z[\Z/n\Z]$ sending $[y]$ to $[y] -
[y+1]$ for any $y \in \Z/n\Z$. Therefore $\Ker b = \Coker b = \Z$
irrespective of $[n]$, and this gives \eqref{4.term} once we
establish functoriality. To do this, consider a map $g:[n] \to [m]$
in $\Lambda$. Then under our identifications $j^o_!\Z([n]) \cong
j_*\Z([n]) \cong \Z[\Z/n\Z]$, $j^o_!\Z([m]) \cong j_*\Z([m]) \cong
\Z[\Z/m\Z]$, the maps $g_!,g_*:\Z[\Z/n\Z] \to \Z[\Z/m\Z]$
corresponding to $j^o_!$ resp. $j_*$ are given by
$$
g_!([y]) = [g(y)], \qquad g_*([y]) = \sum_{y' \in g_\hash^{-1}(y)}
[a'],
$$
and we note that by \eqref{lambda.du}, $g_\hash^{-1}(y)$ is the
(possibly empty) set of residues
$g(y),g(y)+1,\dots,g(y+1)-1$. Therefore we have
$$
b(g_*([y])) = \sum_{y' \in g_\hash^{-1}(y)} ([y'] - [y'-1]) = [g(y)] -
[g(y+1)],
$$
and the right-hand side is exactly $g_!(b([y]))$.
\endproof

\begin{remark}
Lemma~\ref{4-le} has the following geometric meaning. Any $[n] \in
\Lambda$ can be interpreted as a cellular decomposition of the
circle $S^1$ with $0$-cells numbered by elements $y \in
Y([n])$. Morphisms in $\Lambda$ give homotopy classes of cellular
maps, and duality $\Lambda \cong \Lambda^o$ corresponds to taking
the dual cellular decomposition. Then at each $[n] \in \Lambda$, the
middle two terms in \eqref{4.term} give the cellular chain complex
computing the homology $H_\idot(S^1,\Z)$.
\end{remark}

Combining \eqref{E.proj} and \eqref{4.term}, for any $E \in
\Fun(\Lambda,A)$, we obtain a functorial complex
\begin{equation}\label{j.dg}
\begin{CD}
j_*j^*E @>>> j_!^oj^{o*}E
\end{CD}
\end{equation}
of length $2$, with homology objects in degrees $0$ and $1$ 
identified with $E$.

\begin{defn}\label{J.defn}
For any $E \in \Fun(\Lambda,A)$, {\em the complex $\K_\idot(E)$} in
$\Fun(\Lambda,A)$ is the complex \eqref{j.dg}, with $\K_0(E) =
j_!^oj^{o*}E$ and $\K_1(E) = j_*j^*E$.
\end{defn}

For any $l > 1$, the commutative square \eqref{i.l.sq} and the
similar square for the category $\Delta^o$ induce base change
isomorphisms
\begin{equation}\label{bc.i.j}
j^o_{l!} \circ i_l^* \cong i_l^* \circ j^o_!, \qquad i_l^* \circ j_*
\cong j_{l*} \circ i_l^*.
\end{equation}
Thus we can apply $i_l^*$ to the exact sequence \eqref{4.term} and
obtain a four-term exact sequence
$$
\begin{CD}
0 @>>> \Z @>>> j_{l*}\Z @>>> j^o_{l!}\Z @>>> \Z @>>> 0
\end{CD}
$$
and a canonical complex
$$
\begin{CD}
j_{l*}j_l^*E \cong E \otimes j_{l*}\Z @>>> j^o_{l!}j_l^{o*}E \cong E
\otimes j^o_{l!}\Z
\end{CD}
$$
for any $E \in \Fun(\Lambda_l,A)$, with homology objects in degrees
$0$ and $1$ identified with $E$. By abuse of notation, we will
denote this complex by $\K_\idot(E)$, same as in
Definition~\ref{J.defn}. We then have a canonical identification
$$
i_l^*\K_\idot(E) \cong \K_\idot(i_l^*E).
$$
For any object $[n] \in \Lambda$ and any representation $M$ of the
cyclic group $\Aut([n])=\Z/nl\Z$ with generator $\sigma \in
\Z/nl\Z$, we will denote by $\K_\idot(M)$ the complex
\begin{equation}\label{K.M}
\begin{CD}
M @>{\id - \sigma}>> M
\end{CD}
\end{equation}
placed in homological degrees $0$ and $1$. Then for any $E \in
\Fun(\Lambda_l,A)$, we have
\begin{equation}\label{k.n.s}
\K_\idot(E)([n])_\sigma \cong \K_\idot(E([n])).
\end{equation}

\begin{defn}\label{hh.defn}
For any ring $A$, integer $l \geq 1$, and object $E \in
\Fun(\Lambda_l,A)$, the Hochschild resp. cyclic homology groups of
the object $E$ are given by
$$
HH_\idot(E) = H_\idot(\Delta^o,j_l^{o*}E) =
H_\idot(\Lambda_l,j^o_{l!}j^{o*}_lE), \qquad HC_\idot(E) =
H_\idot(\Lambda_l,E).
$$
The {\em Hochschild complex} $CH_\idot(E)$ is the standard chain
complex $C_\idot(j^{o*}_lE)$ of the simplicial $A$-module $j_l^*E$.
\end{defn}

\begin{lemma}\label{la.ho}
For any integer $l \geq 1$, any ring $A$, and any $E \in
\Fun(\Lambda_l,A)$, we have $H_\idot(\Lambda_l,j_{l*}j_l^*E)=0$, so
that the natural map
\begin{equation}\label{flat}
HH_\idot(E) \cong H_\idot(\Lambda_l,j^o_{l!}j_l^{o*}E) \to
HC_\idot(\K_\idot(E))
\end{equation}
is an isomorphism. Moreover, for any $E \in \Fun(\Lambda,A)$, the
adjunction maps
\begin{equation}\label{h.l.i}
HC_\idot(i_l^*E) \to HC_\idot(E), \qquad
HH_\idot(i_l^*E) \to HH_\idot(E)
\end{equation}
are also isomorphisms.
\end{lemma}

\proof{} The first claim is \cite[Lemma 1.10]{ka0}, and the second
claim is \cite[Lemma 1.14]{ka0}.
\endproof

We note that the Hochschild homology isomorphism of \eqref{h.l.i}
holds for any simplicial $A$-module $E$, and it can be realized
explicitly: map \eqref{s.l} induces a map
\begin{equation}\label{h.l}
s_l:\bi_l^*E \to E,
\end{equation}
and this map becomes a quasiisomorphism after we pass to the
standard complexes. In fact, $E$ does not even need to be an
$A$-module --- the map \eqref{h.l} exists equally well for
simplicial objects of any nature, for example for simplicial sets
(it is known as the {\em edgewise subdivision} map and goes back to
\cite{edge}). Realizing explicitly the cyclic homology isomorphism
of \eqref{h.l.i} is more difficult (see e.g.\ \cite[Subsection
  4.1]{coper}).

Since the embedding $j_l:\Delta \to \Lambda_l$ of \eqref{j.l}
extends to an embedding $\wt{j}_l:\Delta^o \times \ppt_l \to
\Lambda_l$, for any $E \in \Fun(\Lambda_l,A)$, the pullback $j_l^*E$
carries a canonical action of the cyclic group $\Z/l\Z$. Then by
adjunction, $j_{l*}\Z$ also carries a $\Z/l\Z$-action. Analogously,
we have a natural $\Z/l\Z$-action on the functors $j^{o*}_l$ and on
$j^o_{l!}\Z$. This is compatible with the differential and turns
$\K_\idot(E) \cong E \otimes K_\idot(\Z)$ into a complex of
$\Z/l\Z$-modules in $\Fun(\Lambda_l,A)$.

\begin{corr}\label{i.l.triv}
For any $E \in \Fun(\Lambda_l,A)$, the natural $\Z/l\Z$-action on
the Hochschild homology groups $HH_\idot(E)$ induced by the
$\Z/l\Z$-action on $j^{o*}_lE$ is trivial.
\end{corr}

\proof{} Let $B_0 = \Z[\Z/l\Z]$ be the group algebra of the cyclic
group $\Z/l\Z$, with $[\sigma] \in B_0$ being the class of the
generator $\sigma \in \Z/l\Z$, and let $B_\idot = B_0(\eps)$ be the
free commutative DG algebra over $B_0$ generated by one generator
$\eps$ of homological degree $1$, with $\eps^2 = 0$ and $d\eps = 1 -
[\sigma]$. Then in terms of the identification \eqref{flat}, the
$\Z/l\Z$-action on $HH_\idot(E)$ is induced by the $B_0$-module
structure on the complex $\K_\idot(E) \cong E \otimes \K_\idot(\Z)$,
and it suffices to prove that $\K_\idot(\Z)$ is actually a DG module
over the whole $B_\idot$. To do this, we have to extend the action
map $B_0 \otimes \K_0(\Z) \to \K_0(\Z)$ to a map of complexes
\begin{equation}\label{B.act}
B_\idot \otimes K_0(\Z) \to \K_\idot(\Z)
\end{equation}
in such a way that the resulting $B_\idot$-action is associative and
compatible with the differentials. However, by adjunction, we have
an identification
$$
\begin{aligned}
\Hom(V \otimes \K_0(\Z),F) &= \Hom(V \otimes j^o_{l!}\Z,F) \cong\\
&\cong \Hom(V \otimes \Z,j^{o*}_lF) \cong \Hom(V,F([1]))
\end{aligned}
$$
for any abelian group $V$ and any $F \in \Fun(\Lambda_l,\Z)$, so that
it suffices to construct the extended map \eqref{B.act} and check
that it defines a $B_\idot$-module structure on
$\K_\idot(\Z)([1])$. This is easy, since in fact $\K_\idot(\Z)([1])
\cong B_\idot$.
\endproof

Now for any ring $A$, integer $l \geq 1$, and $E \in
\Fun(\Lambda_l,A)$, consider the adjunction maps
\begin{equation}\label{kappa}
\kappa_0:\K_0(E) = j^o_{l!}j^{o*}_lE \to E, \qquad
\kappa_1:E \to j^o_{l*}j^{o*}_lE = \K_1(E),
\end{equation}
and let
\begin{equation}\label{B.0}
B = \kappa_1 \circ \kappa_0:\K_0(E) \to E \to \K_1(E)
\end{equation}
be their composition. We can then form a natural periodic resolution
\begin{equation}\label{reso}
\begin{CD}
@>{B}>> \K_1(E) @>>> \K_0(E) @>{B}>> \K_1(E) @>>> \K_0(E)
\end{CD}
\end{equation}
of the object $E$. If we denote by $u$ the endomorphism of
\eqref{reso} obtained by shifting to the left by $2$ terms, then $u$
gives a natural class $u \in \Ext^2(E,E)$ (equivalently, this is the
class represented by Yoneda by the complex $\K_\idot(E)$). Taking
the cone of $u$ and using the isomorphism \eqref{flat}, we obtain
the {\em Connes long exact sequence}
\begin{equation}\label{con.seq}
\begin{CD}
HH_\idot(E) @>>> HC_\idot(E) @>{u}>> HC_{\idot-2}(E) @>>>
\end{CD}
\end{equation}
On the other hand, taking the odd-numbered terms of the stupid
filtration on the complex \eqref{reso}, we obtain a natural spectral
sequence
\begin{equation}\label{sp.se}
HH_\idot(E)[u^{-1}] \Rightarrow HC_\idot(\Lambda_l,E),
\end{equation}
where the left-hand side is shorthand for ``formal polynomials in
one variable $u^{-1}$ of homological degree $2$ with coefficients in
$HH_\idot(E)$''. This is known as the {\em Hochschild-to-cyclic}, or
{\em Hodge-to-de Rham} spectral sequence.  The first non-trivial
differential
\begin{equation}\label{c.ts}
B:HH_\idot(E) \to HH_{\idot+1}(E)
\end{equation}
is known as the {\em Connes-Tsygan differential}, or sometimes {\em
  Rinehart differential}. In terms of the identification
\eqref{flat}, this differential is induced by the map
\begin{equation}\label{B.gen}
B:\K_\idot(E) \to \K_\idot(E)[-1]
\end{equation}
which is in turn induced by the natural map $B$ of \eqref{B.0}
(this is why we use the same notation for all three).

\section{Preliminaries II.}\label{prem2.sec}

\subsection{Trace maps.}

Let us now explore the homological properties of the projections
$\pi_\idot$ of \eqref{i.pi.l}. Fix an integer $l \geq 1$, and
consider the projection $\pi_l:\Lambda_l \to \Lambda$. This is a
bifibration with fiber $\ppt_l$, so that we have functorial maps
\eqref{tr.gamma} and \eqref{e.gamma}. We simplify notation by
writing
\begin{equation}\label{e.tr.pi}
\e_l = \e_{\pi_l}:\pi_{l*}E \to \pi_{l!}E, \qquad
\tr_l=\tr_{\pi_l}:\pi_{l!}E \to \pi_{l*}E.
\end{equation}
If $E \in \Fun(\Lambda_l,A)$ is $\pi_l$-free, then the map $\tr_l$
of \eqref{e.tr.pi} is an isomorphism, and $E$ is acyclic both for
the left-exact functor $\pi_{l!}$ and for the right-exact functor
$\pi_{l*}$ (that is, $L^i\pi_{l!}E = \pi_{l!}E$, $R^i\pi_{l*}E =
\pi_{l*}E$). In particular, by \eqref{k.exp} and \eqref{bc.i.j}, this
applies to functors of the form $j^o_{l!}(E)$, $E \in
\Fun(\Delta^o,A)$ and $j_{l*}(E)$, $E \in \Fun(\Delta,A)$, so that
we have a natural isomorphism
\begin{equation}\label{tr.k}
\tr_l:\pi_{l!}\K_\idot(E) \cong \pi_{l*}\K_\idot(E)
\end{equation}
for any $E \in \Fun(\Lambda_p,E)$. Moreover, since
$\pi_{l!}\K_\idot(E) \cong L^\hdot\pi_{l!}\K_\idot(E)$, we have
natural identifications
\begin{equation}\label{pi.k}
HC_\idot(\K_\idot(E)) \cong HC_\idot(\pi_{l!}\K_\idot(E)) \cong
HC_\idot(\pi_{l*}\K_\idot(E)).
\end{equation}
Under the isomorphism \eqref{flat}, all these groups are further
identified with $HH_\idot(E)$. Moreover, to simplify notation, we
introduce the following.

\begin{defn}\label{K.l.defn}
For any $E \in \Fun(\Lambda,A)$ and any $l \geq 1$, the complex
$\K^l_\idot(E)$ in $\Fun(\Lambda,A)$ is given by
$$
\K^l_\idot(E) = \pi_{l!}\K_\idot(i_l^*E) \cong
\pi_{l*}\K_\idot(i_l^*E).
$$
\end{defn}

Then by definition, for any $E \in \Fun(\Lambda,A)$ and any $l \geq
1$, the complex $\K^l_\idot(E)$ fits into an exact sequence
\begin{equation}\label{K.l.4}
\begin{CD}
0 @>>> \pi_{l*}E @>>> \K_1^l(E) @>>> \K^l_0(E)
@>>> \pi_{l!}E @>>> 0,
\end{CD}
\end{equation}
and \eqref{pi.k} together with the isomorphism of Lemma~\ref{la.ho}
provides a canonical isomorphism
\begin{equation}\label{K.l.cor}
HH_\idot(E) \cong HC_\idot(\K^l_\idot(E)).
\end{equation}

\begin{lemma}\label{k.l.l}
For any $l \geq 1$ and any $E \in \Fun(\Lambda,A)$, we have
natural isomorphism of complexes
$$
\K_\idot(E) \cong \pi_{l!}\K_\idot(\pi_l^*E) \cong
\pi_{l*}\K_\idot(\pi_l^*E).
$$
\end{lemma}

\proof{} By the projection formula \eqref{proj.eq}, we have
$\pi_{l!}\K_\idot(\pi_l^*E) \cong E \otimes \pi_{l!}\pi_l^*\Z$ and
$\pi_{l*}\K_\idot(\pi_l^*E) \cong E \otimes \pi_{l*}\pi_l^*\Z$, so
that it suffices to consider the case $E = \Z$. In this case, the
claim immediately follows from \eqref{k.exp} and \eqref{bc.i.j}.
\endproof

In particular, for any $E \in \Fun(\Lambda_l,A)$, the adjunction
maps $E \to \pi_l^*\pi_{l!}E$, $\pi_l^*\pi_{l*}E \to E$ induce
natural maps
\begin{equation}\label{nu.phi.eq}
\begin{aligned}
\nu_l:&\pi_{l!}\K_\idot(E) \to \pi_{l!}\K_\idot(\pi_l^*\pi_{l!}E) \cong
\K_\idot(\pi_{l!}E),\\
\phi_l:&\K_\idot(\pi_{l*}E) \to \pi_{l*}\K(E).
\end{aligned}
\end{equation}

\begin{lemma}\label{nu.phi}
For any $l \geq 1$, ring $A$, and $E \in \Fun(\Lambda_l,A)$, the
map $\nu_l$ is an isomorphism on homology in degree $0$ and equal to
the map $\e_l$ of \eqref{e.tr.pi} in degree $1$, and $\phi_l$ is an
isomorphism on homology in degree $1$, and equal to $\e_l$ in degree
$0$.
\end{lemma}

\proof{} The canonical maps $\kappa_0$ of \eqref{kappa} induce
isomorphisms of the $0$-th homology objects of the complexes
$\K_\idot(E)$ resp. $\K_\idot(\pi_{l!}E)$ with $E$
resp. $\pi_{l!}E$. Since the functor $\pi_{l!}$ is right-exact,
$\pi_{l!}(\kappa_0)$ then gives an isomorphism between $0$-th
homology of $\pi_{l!}\K_\idot(E)$ and $\pi_{l!}E$, and since by
definition, we have $\pi_{l!}(\kappa_0) = \kappa_0 \circ \nu_l$,
$\nu_l$ induces an isomorphism on homology in degree
$0$. Analogously, $\phi_l$ induces an isomorphism on homology in
degree $1$. Then by \eqref{tr.k}, the complexes
$\pi_{l!}\K_\idot(E)$ and $\pi_{l*}\K_\idot(E)$ are one and the
same, so that to prove that $\nu_l$ induces the map $\e_l$ on
homology in degree $1$, we have to prove that the diagram
$$
\begin{CD}
\pi_{l*}E @>{\pi_{l*}(\kappa_1)}>> \pi_{l*}\K_1(E) \cong
\pi_{l!}\K_1(E)\\
@V{\e_l}VV @VV{\nu_l}V\\
\pi_{l!}E @>{\kappa_1}>> \K_1(\pi_{l!}E)
\end{CD}
$$
is commutative. But $\pi_{l*}\K_1(E) = \pi_{l*}j_{l*}j_l^*E
\cong j_*j_l^*E$ and $\K_1(\pi_{l!}E) = j_*j^*\pi_{l!}E$, so that by
adjunction, it is equivalent to proving that the diagram
$$
\begin{CD}
j^*\pi_{l*}E @>>> j_l^*E\\
@V{\e_l}VV @VVV\\
j^*\pi_{l!}E @= j^*\pi_{l!}E
\end{CD}
$$
is commutative -- in other words, that the composition $j^*\pi_{l*}E
\to j_l^*E \to j^*\pi_{l!}E$ is the map $j^*(\e_l)$. After evaluating
at an object $[n] \in \Delta$, this composition reads as
$$
\begin{CD}
E([n])^\sigma @>>> E([n]) @>>> E([n])_\sigma,
\end{CD}
$$
and it is equal to $\e_l$ by definition. The argument
for $\phi_l$ is dual.
\endproof

We will also need a result on compatibility between the
Connes-Tsygan differential \eqref{c.ts} and the trace map
\eqref{e.tr.pi}. Consider an object $E$ in the category
$\Fun(\Lambda_l,A)$, and denote by $b:\K_\idot(\pi_{l!}E) \to
\K_\idot(\pi_{l*}E)[-1]$ the composition map
\begin{equation}\label{b.B.eq}
\begin{CD}
\K_\idot(\pi_{l!}E) @>{\kappa_0}>> \pi_{l!}E @>{\tr_l}>> \pi_{l*}E
@>{\kappa_1}>> \K_\idot(\pi_{l*}E)[-1],
\end{CD}
\end{equation}
where $\kappa_0$, $\kappa_1$ are the adjunction maps
\eqref{kappa}.

\begin{lemma}\label{fdv.loc}
Under the identifications \eqref{pi.k}, the composition map
$$
\begin{CD}
\pi_{l!}\K_\idot(E) @>{\nu_l}>> \K_\idot(\pi_{l!}E) @>{b}>>
\K_\idot(\pi_{l*}E)[-1] @>{\phi_l}>> \pi_{l*}\K_\idot(E)[-1]
\end{CD}
$$
induces the Connes-Tsygan differential $B$ of \eqref{c.ts}.
\end{lemma}

\proof{} The Connes-Tsygan differential is by definition induced by
the composition
$$
\begin{CD}
\pi_{l!}\K_\idot(E) @>{\pi_{l!}(B)}>> \pi_{l!}\K_\idot(E)[-1] @>{\tr_l}>>
\pi_{l*}\K_\idot(E)[-1],
\end{CD}
$$
where $B$ is as in \eqref{B.gen}, and $\tr_l$ on the right is the
canonical isomorphism \eqref{tr.k}. By \eqref{B.0}, this composition
is in turn given by the composition
$$
\begin{CD}
\pi_{l!}\K_0(E) @>{\pi_{l!}(\kappa_0)}>> \pi_{l!}E
@>{\pi_{l!}(\kappa_1)}>> \pi_{l!}\K_1(E) @>{\tr_l}>> \pi_{l*}\K_1(E).
\end{CD}
$$
Since the map $\tr_l$ is functorial, this is equal to the
composition
$$
\begin{CD}
\pi_{l!}\K_0(E) @>{\pi_{l!}(\kappa_0)}>> \pi_{l!}E @>{\tr_l}>> \pi_{l*}E
@>{\pi_{l*}(\kappa_1)}>> \pi_{l*}\K_1(E),
\end{CD}
$$
and by definition, $\pi_{l!}(\kappa_0) = \kappa_0 \circ \nu_l$ and
$\pi_l^*(\kappa_1) = \phi_l \circ \kappa_1$.
\endproof

\subsection{Extra structures.}\label{mult.subs}

If the ring $A$ is commutative, then the categories
$\Fun(\Lambda,A)$, $\Fun(\Delta^o,A)$ acquire natural pointwise
tensor products, and the derived categories $\D(\Lambda,A)$,
$\D(\Delta^o,A)$ in turn acquire the derived tensor product
$\lotimes$. It is well-known that the Hochschild homology functor is
multiplicative: for any two objects $E,E' \in \D(\Lambda,A)$, we
have a natural K\"unneth quasiisomorphism
\begin{equation}\label{kunn}
CH_\idot(E) \lotimes CH_\idot(E') \cong CH_\idot(E \lotimes E'),
\end{equation}
where we denote by $CH_\idot(-)$ the object in the derived category
$\D(A\amod)$ corresponding to the Hochschild homology $HH_\idot(-)$
(to make things more canonical, one can agree to take as
$CH_\idot(E)$ the standard complex of the simplicial $A$-module
$j^{o*}E$). The Connes-Tsygan differential $B$ of \eqref{c.ts} is
compatible with the K\"unneth quasiisomorphism \eqref{kunn} -- that
is, we have
\begin{equation}\label{B.kun}
B_{E \lotimes E'} = B_{E} \otimes \id + \id \otimes B_{E'}.
\end{equation}
For a characteristic-independent proof of this, see
e.g.\ \cite[Lemma 1.3]{cart}.

We will also need to consider one additional structure on cyclic
objects essentially introduced in \cite{witt}. Fix a prime $p$, and
consider the functors $i_p,\pi_p:\Lambda_p \to \Lambda$ and the map
$\tr_p$ of \eqref{e.tr.pi}. Similarly, denote
$\tr^\dg_p=\tr^\dg_{\pi_p}$, where $\tr^\dg$ is the map
\eqref{tr.dg}

\begin{defn}\label{FV.def}
An {\em $FV$-structure} on an object $E \in \Fun(\Lambda,A)$ is
given by two maps $V:i^*_pE \to \pi_p^*E$, $F:\pi_p^*E \to i_p^*E$
such that
\begin{equation}\label{FV.tr}
F \circ V = \tr^\dg_p:i_p^*E \to i_p^*E.
\end{equation}
\end{defn}

\begin{remark}
Formally, Definition~\ref{FV.def} makes sense for any $p \geq 1$;
however, the notion is meaningful only when $p$ is a prime.
\end{remark}

For any object $E \in \Fun(\Lambda,A)$ equipped with an
$FV$-structure $\langle V,F \rangle$, the maps $V$ and $F$ induce by
adjunction natural maps
\begin{equation}\label{bar.FV}
\overline{V}:\pi_{p!}i_p^*E \to E, \qquad \overline{F}:E \to
\pi_{p*}i_p^*E
\end{equation}
whose composition is equal to $\tr_p$. We can consider the
compositions
\begin{equation}\label{compo}
\begin{CD}
\pi_{p!}\K_\idot(i_p^*E) @>{\nu_p}>> \K_\idot(\pi_{p!}E)
@>{\K_\idot(\overline{V})}>> \K_\idot(E),\\
\K_\idot(E) @>{\K_\idot(\overline{F})}>> \K_\idot(\pi_{p*}E) @>{\phi_p}>>
\pi_{p*}\K_\idot(i_p^*E),
\end{CD}
\end{equation}
where $\nu_p$, $\phi_p$ are the maps \eqref{nu.phi.eq}, and by
Lemma~\ref{la.ho} together with the isomorphism \eqref{K.l.cor},
these induce natural maps
\begin{equation}\label{FV.HH}
V,F:HH_\idot(E) \to HH_\idot(E).
\end{equation}

\begin{lemma}\label{FdV.le}
For any $E \in \Fun(\Lambda,A)$ equipped with an $FV$-structure, we have
$$
FV=p\id:HH_\idot(E) \to HH_\idot(E)
$$
and
$$
FBV = B:HH_\idot(E) \to HH_{\idot+1}(E),
$$
where $V$ and $F$ are the maps \eqref{FV.HH}, and $B$ is the
differential \eqref{c.ts}.
\end{lemma}

\proof{} The first identity immediately follows from \eqref{FV.tr}
and Corollary~\ref{i.l.triv}. To prove the second equality, note
that by Lemma~\ref{fdv.loc} and \eqref{compo}, it suffices to prove
that the canonical map $b:\K_0(\pi_{p!}i_p^*E) \to
\K_1(\pi_{p*}i_p^*E)$ given by \eqref{b.B.eq} coincides with the
composition
$$
\begin{CD}
\K_0(\pi_{p!}i_p^*E) @>{\K_0(V^\dg)}>> \K_0(E) @>{B}>> \K_1(E)
@>{\K_1(F^\dg)}>> \K_1(\pi_{p*}i_p^*E).
\end{CD}
$$
By definition, we have $B = \kappa_1 \circ \kappa_0$, and we also
have $\tr_p=F^\dg \circ V^\dg$. Thus it suffices to prove that the
diagrams
$$
\begin{CD}
\K_0(\pi_{p!}i_p^*E) @>{\kappa_0}>>
\pi_{p!}i_p^*E\\
@V{\K_0(V^\dg)}VV @VV{V^\dg}V\\
\K_0(E) @>{\kappa_0}>> E
\end{CD}
\qquad
\begin{CD}
E @>{\kappa_1}>> \K_1(E)\\
@V{F^\dg}VV @VV{\K_1(F^\dg)}V\\
\pi_{p*}i_p^*E @>{\kappa_1}>>
\K_1(\pi_{p*}i_p^*E)
\end{CD}
$$
are commutative. This immediately follows from the functoriality of
the adjunction maps $\kappa_0$, $\kappa_1$.
\endproof

We note that for any $E \in \Fun(\Lambda,A)$ equipped with an
$FV$-structure $\langle V,F \rangle$, we can iterate the maps $V$
and $F$. Namely, for any integer $m \geq 1$, we can inductively
apply base change to the Cartesian square \eqref{la.ml} and obtain
natural maps
\begin{equation}\label{FV.ite}
V^m:i_{p^m}^*E \to \pi_{p^m}^*E, \qquad F^m:\pi_{p^m}^*E \to
i_{p^m}^*E.
\end{equation}
By adjunction, these induce maps
\begin{equation}\label{bar.FV.ite}
\overline{V}^m:\pi_{p^m!}i_{p^m}^*E \to E, \qquad \overline{F}^m:E \to
\pi_{p^m*}i_{p^m}^*E,
\end{equation}
a generalization of \eqref{bar.FV}, and taking the compositions
\eqref{compo}, we obtain natural maps
$$
V^m,F^m:HH_\idot(E) \to HH_\idot(E).
$$
By base change, these are indeed iterates of the maps \eqref{FV.HH},
so that our notation is consistent.

\subsection{Yoneda sets.}

Recall now that sending an object $[n] \in \Lambda$ to the set of
objects of the category $[n]_\Lambda$ defines a functor $Y:\Lambda
\to \Sets$, and this functor is corepresented by $[1] \in \Lambda$. One
can also consider the functors corepresented by different objects
$[n] \in \Lambda$; we now prove some results on their structure that
we will need in Subsection~\ref{free.subs}.

For any $[n] \in \Lambda$, with automorphism group $\Aut([n]) =
\Z/n\Z$, denote by $e^n:\ppt_{n} \to \Lambda$ the embedding onto
$[n]$. Consider the functor
\begin{equation}\label{e.n}
e^n_!:A[\Z/n\Z]\amod \cong \Fun(\ppt_n,A) \to
\Fun(\Lambda,A).
\end{equation}
To describe this functor explicitly, let $\Lambda^n$ be the category
of maps $f:[n] \to [m]$ in $\Lambda$, with maps from $f':[n] \to
[m']$ to $f:[n] \to [m]$ given by commutative diagrams
$$
\begin{CD}
[n] @>{f'}>> [m']\\
@V{\sim}VV @VVV\\
[n] @>{f}>> [m].
\end{CD}
$$
Then sending $f:[n] \to [m]$ to $[n]$ resp.\ $[m]$ gives functors
$x:\Lambda^n \to \ppt_n$, $y_n:\Lambda^n \to \Lambda$, the embedding
$e^n$ factors as $e^n = y_n \circ \wt{e}^n$ with $\wt{e}^n$
left-adjoint to $x$, and we have
$$
e^n_! \cong y_{n!} \circ \wt{e}^n_! \cong y_{n!} \circ x^*.
$$
The functor $y^n:\Lambda^n \to \Lambda$ is a discrete
cofibration. Explicitly, define a functor $\wt{Y}^n:\Lambda \to
\Sets$ by
$$
\wt{Y}^n([m]) = \Lambda([n],[m]), \qquad [m] \in \Lambda,
$$
and let $Y^n=\wt{Y}^n/\Aut([n])$. Then the $\Aut([n])$-action on
$\wt{Y}^n$ is free, and $Y^n$ corresponds by the Grothendieck
construction to the discrete cofibration $y^n$.  If $n=1$, then
$\wt{Y}^1=Y^1$ coincides with the functor $Y:\Lambda \to \Sets$, the
category $\Lambda^1$ is equivalent to $\Delta^o$, and $y^1:\Lambda^1
\to \Lambda$ is then identified with $j^o:\Delta^o \to \Lambda$.

More generally, take some integer $l \geq 1$, denote by
$e^n_l:\ppt_{nl} \to \Lambda_l$ the embedding onto $[n] \in
\Lambda_l$, and consider the functor
\begin{equation}\label{e.n.l}
e^n_{l!}:A[\Z/nl\Z]\amod \cong \Fun(\ppt_{nl},A) \to
\Fun(\Lambda_l,A)
\end{equation}
that reduces to \eqref{e.n} when $l=1$. Then as in the case $l=1$,
we can define the functor $Y^n_l:\Lambda_l \to \Sets$ by
$$
Y^n_l([m]) = \Lambda_l([n],[m])/\Aut([n]), \qquad [m] \in \Lambda_l,
$$
consider the corresponding discete cofibration $y^n_l:\Lambda^n_l
\to \Lambda_l$, and note that $e^n_l \cong y^n_l \circ \wt{e}^n_l$ for
some embedding $\wt{e}^n_l:\ppt_{nl} \to \Lambda^n_l$ that admits a
right-adjoint projection $x:\Lambda^n_l \to \ppt_{nl}$. We then
again have an isomorphism
\begin{equation}\label{e.x}
e^n_{l!} \cong y^n_{l!} \circ x^*.
\end{equation}

\begin{lemma}\label{e.n.le}
For any $[n] \in \Lambda_l$, the functor \eqref{e.n.l} is
exact. Moreover, for any $A[\Z/nl\Z]$-module $M$, the Hochschild
complex $CH_\idot(e^n_{l!}M)$ is naturally quaisisomorphic to the
complex $\K_\idot(M)$ of \eqref{K.M}, and the Connes-Tsygan
differential $B$ is induced by the trace map $\tr_{\Z/nl\Z}$.
\end{lemma}

\proof{} Since $y^n_l:\Lambda^n_l \to \Lambda_l$ is a discrete
cofibration, the first claim immediately follows from
\eqref{e.x}. Moreover, we tautologically have
\begin{equation}\label{e.n.l.H}
HC_\idot(e^n_{l!}M) = H_\idot(\Lambda_l,e^n_{l!}E) \cong
H_\idot(\Z/nl\Z,M),
\end{equation}
and by \eqref{k.n.s}, if we compute $H_\idot(\Z/nl\Z,M)$ by the
standard $2$-periodic resolution, then the isomorphism
\eqref{e.n.l.H} is compatible with the periodicity. Then the second
claim immediately follows from \eqref{con.seq}.
\endproof

Now note that for any $n,l \geq 1$, we have $i_l \circ e^n_l \cong
e^{nl}$, so that by adjunction, we have a functorial map
\begin{equation}\label{i.n.l}
e^n_{l!} \to i_l^* \circ i_{l!} \circ e^{nl}_! \cong i_l^* \circ
e^{nl}_!.
\end{equation}
Then for any $A[\Z/nl\Z]$-module $M$, the induced map
$$
HH_\idot(e^n_{l!}M) \to H_\idot(i_l^*e^{nl}_!M)
$$
is an isomorphism by Lemma~\ref{la.ho}. What we will need is the
following refinement of this fact.

\begin{prop}\label{Y.prop}
For any $n,l \geq 1$ and $M \in A[\Z/nl\Z]\amod$, the map
$$
HH_\idot(\pi_{l!}e^n_{l!}M) \to HH_\idot(\pi_{l!}i_l^*e^{nl}_!M)
$$
induced by the map \eqref{i.n.l} is an isomorphism.
\end{prop}

In order to prove this, we need to analyse the structure of the
functors $Y^n_l$ in some detail. For any integer $m \geq 1$, denote
\begin{equation}\label{Y.nm}
Y^{n,m}_l = Y^n_l \circ i_m:\Lambda_{ml} \to \Sets,
\end{equation}
and let $y^{n,m}_l:\Lambda^{n,m}_l \to \Lambda_{ml}$ be the
corresponding discrete cofibration. To simplify notation, set
$Y^{n,m}=Y^{n,m}_1$, $\Lambda^{n,m}=\Lambda^{n,m}_1$,
$y^{n,m}=y^{n,m}_1$. Then for any $n,m,l \geq 1$, the embeddings
\eqref{i.l.map} taken together define a canonical embedding
\begin{equation}\label{i.l.Y}
\phi^{n,m}_l:Y^{n,m}_l \subset Y^{nl,ml}
\end{equation}
and a functor $\phi^{n,m}_l:\Lambda^{n,m}_l \to \Lambda^{nl,ml}$
such that $y^{nl,ml} \circ \phi^{n,m}_l \cong y^{n,m}_l$. We have a
natural projection $x:\Lambda^{n,m}_l \to \ppt_{nl}$ such that
$$
i_m^* \circ e^n_{l!} \cong y^{n,m}_{l!} \circ x^*,
$$
where $i_m:\Lambda_l \to \Lambda_{ml}$ is the functor \eqref{la.ml},
and in terms of this identification, the map \eqref{i.n.l} is
induced by the embedding $\phi^{n,1}_l:\Lambda^n_l \to
\Lambda^{nl,l}$.

\proof[Proof of Proposition~\ref{Y.prop}.] Denote
$E=i_l^*e^{nl}_!M$, and for any $r \geq 1$ dividing $l$, denote $E_r
= i_m^*e^{nm,m}_{r!}M$, where $m=l/r$. Note that the embeddings
\eqref{i.l.Y} induce injective maps $\nu_r:E_r \to E$. Filter $E$ by
setting
\begin{equation}\label{flt}
F_iE = \sum_{r|l,l/r \leq i}\nu_r(E_r) \subset E.
\end{equation}
By definition, we have $F_lE=E$, and $F_1E = E_l = e^n_{l!}M$, with
$\nu_l:E_l \to E$ being the map \eqref{i.n.l}. We have $\gr_iFE=0$
unless $i$ divides $l$, and for any $m|l$, $r=l/m$, we have $\gr^m_F
E \cong \gr^m_F E_r$, where $E_r \subset E$ is equipped with the
induced filtration. For any $r|l$, the map $HH_\idot(E_l) \to
HH_\idot(E_r)$ induced by the embedding $E_l \subset E_r$ is an
isomorphism by Lemma~\ref{la.ho}, and by induction on $r$, this
implies that
\begin{equation}\label{hh.gr.m}
HH_\idot(\gr^m_FE)=0
\end{equation}
for any $m|l$, $m > 1$.

Now note that the filtration \eqref{flt} is induced by the
filtration $F_\idot$ on $Y^{nl,l}$ given by the unions of the images
of the embeddings \eqref{i.l.Y}, and after evaluating at any object
$[s] \in \Lambda_l$, the latter filtration splits: we have
\begin{equation}\label{splt}
Y^{nl,l}([s]) = \coprod_{m \leq i,m |
  l}\overline{Y}^{nm,m}_r([s]),
\end{equation}
where we denote $\overline{Y}^{nm,m}_r([s]) = Y^{nm,m}_r([s])
\setminus (F_{m-1}Y^{nl,l}([s]) \cap Y^{nm,m}_r([s]))$. Moreover, by
Lemma~\ref{fp.le}, the splitting \eqref{splt} is invariant under the
action of $\Z/l\Z \subset \Aut([s])$. Therefore the filtration
$F_\idot E([s])$ admits a $\Z/l\Z$-equivariant splitting. By base
change, this implies that the map $\pi_{l!}F_iE \to \pi_{l!}E$ is
injective for any $i$, and $\pi_{l!}E$ acqures a filtration
$F_\idot$ such that
$$
\gr^i_F\pi_{l!}E \cong \pi_{l!}\gr^i_FE, \qquad 1 \leq i \leq l.
$$
Then to prove the Proposition, it suffices to prove that
$HH_\idot(\pi_{l!}\gr^i_FE)=0$ for any $i > 1$. Since this is
tautologically true if $\gr^i_FE=0$, we may further assume that
$i=m$ divides $l$.

We now observe that by Lemma~\ref{fp.le}, we have a commutative
diagram with Cartesian squares
$$
\begin{CD}
\ppt_{nl} @<{x}<< \Lambda^{nm,m}_r @>{y^{nm,m}_r}>> \Lambda_l\\
@V{\pi_r}VV @VVV @VV{\pi_r}V\\
\ppt_{ml} @<{x}<< \Lambda^{nm,m} @>{y^{nm,m}}>> \Lambda_m,
\end{CD}
$$
where the rightmost vertical arrow is the functor \eqref{la.ml}, the
leftmost arrow denotes by abuse of notation the tautological
projection $\pi_r:\ppt_{nl} \to \ppt_{nm}$, and the arrow in the
middle is a bifibration with fiber $\ppt_r$. Therefore if we
denote $E'=e^n_{m!}\pi_{r*}M$, then by base change, we have
$$
\pi_{r!}E_r = \pi_{r!}y^{nm,m}_{r!}x^*M \cong y^{nm,n}_!x^*\pi_{r*}M
= E',
$$
so that
$$
\pi_{l!}\gr^m_FE \cong \pi_{m!}\pi_{r!}\gr^m_FE \cong
\pi_{m!}\pi_{r!}\gr^m_FE_r \cong \pi_{m!}\gr^m_F\pi_{r!}E_r \cong
\pi_{m!}\gr^m_FE',
$$
where we define the filtration $F_\idot$ on $E'$ in the same way as
on $E$. Thus by induction, we may further assume that $l=m$.

It remains to notice that by Lemma~\ref{fp.le}, for any $[s] \in
\Lambda_m$, the action of $\Z/l\Z \subset \Aut([s])$ on
$\overline{Y}^{nl,l}([s])$ is free. Therefore the object $\gr^l_FE
\in \Fun(\Lambda_l,A)$ is $\pi_l$-free, and
$$
HH_\idot(\pi_{l!}\gr^l_FE) \cong HH_\idot(L^\hdot\pi_{l!}\gr^l_FE)
\cong H_\idot(\Delta^o \times \ppt_l,\gr^l_FE).
$$
This vanishes by \eqref{hh.gr.m}.
\endproof

\section{Constructions.}\label{cons.sec}

\subsection{Algebras.}\label{alg.subs}

Fix a commutative ring $k$. To any associative unital algebra $A$
over $k$ one associates a canonical object $A^{\hush} \in
\Fun(\Lambda,k)$ as follows:
\begin{itemize}
\item on objects, $A^{\hush}([n]) = A^{\otimes_k Y([n])} = A^{\otimes_k
  n}$, with copies of $A$ numbered by elements $y \in Y([n])$,
\item for any map $f:[n'] \to [n]$, the map
$$
A^{\hush}(f):A^{\otimes_k n'} = \bigotimes_{y \in
  Y([n])}A^{\otimes_k f^{-1}(y)} \to A^{\otimes n}
$$
is given by
\begin{equation}\label{a.hash.maps}
A^{\hush}(f) = \bigotimes_{y \in Y([n])}m_{f^{-1}(y)},
\end{equation}
where $m_{f^{-1}(y)}:A^{\otimes_k f^{-1}(y)} \to A$ is the map which
multiplies the entries in the canonical order.
\end{itemize}
Note that \eqref{a.hash.maps} makes sense in an arbitrary symmetric
monoidal category. In particular, for any associative monoid $G$, we
can define the natural functor $G^\hush:\Lambda \to \Sets$ by
setting $G^\hush([n])=G^n$, and with the structure maps given by
\eqref{a.hash.maps}. We will say that a $k$-algebra $A$ is {\em
  monomial} if $A=k[G]$ for such a monoid $G$. In this case, we have
\begin{equation}\label{A.G}
A^\hush \cong k[G^\hush].
\end{equation}
Now take a general associative $k$-algebra $A$, and assume that it
is flat as a $k$-module. Then by definition, the {\em Hochschild
  homology} $HH_\idot(A)$ and the {\em cyclic homology}
$HC_\idot(A)$ are given by
$$
HH_\idot(A) = HH_\idot(A^\hush), \qquad HC_\idot(A) =
HC_\idot(A^{\hush}).
$$
The spectral sequence \eqref{sp.se} reads as
\begin{equation}\label{hh.hc}
HH_\idot(A)[u] \Rightarrow HC_\idot(A).
\end{equation}
Note that by Lemma~\ref{la.ho}, one can equally well replace
$\Lambda$ with $\Lambda_l$, $l \geq 1$, and $A_{\hush}$ with
$i_l^*A_{\hush}$; the resulting spectral sequence is the same.

For any two flat $k$-algebras $A$, $B$, we have a natural
isomorphism
\begin{equation}\label{kun.loc}
(A \otimes_k B)^\hush \cong A^\hush \otimes B^\hush,
\end{equation}
and the K\"unneth quasiisomorphism provides natural multiplication
maps
\begin{equation}\label{hh.mult}
HH_\idot(A) \otimes_k HH_\idot(B) \to HH_\idot(A \otimes_k B).
\end{equation}
By \eqref{B.kun}, the Connes-Tsygan differential $B$ is a derivation
with the respect to these multiplication maps.

In general, both $HH_\idot(A)$ and $HC_\idot(A)$ are just
$k$-modules. However, if $A$ is commutative, then the product map
$m:A \otimes_k A \to A$ is a map of algebras. Then \eqref{kun.loc}
together with the map $m^\hush:(A \otimes_k A)^\hush \to A^\hush$
turns $A^\hush$ into an associative commutative algebra object in
$\Fun(\Lambda,k)$, and the multiplication maps \eqref{hh.mult} turn
$HH_\idot(A)$ into an an associative commutative $k$-algebra. The
Connes-Tsygan differential $B$ is then a derivation with respect to
this algebra structure.

The main comparison theorem for Hochschild homology in the
commutative case is the classic result of Hochschild, Kostant, and
Rosenberg.

\begin{theorem}[\cite{HKR}]\label{hkr}
Assume that the algebra $A$ is commutative and finitely generated
over $k$, and that $X = \Spec A$ is smooth over $k$. Then we have
natural isomorphisms
\begin{equation}\label{hkr.iso}
HH_i(A) \cong \Omega^i_A = H^0(X,\Omega_X^i)
\end{equation}
for any $i \geq 0$, and these isomorphisms are compatible with
multiplication.\endproof
\end{theorem}

We want to emphasize that the Hochschild-Kostant-Rosenberg Theorem
requires no assumptions on $\cchar k$. As for cyclic homology, the
main result is the following (for a characteristic-independent
proof, see e.g.\ \cite[Theorem 2.2]{cart}, or an essentially
equivalent \cite[Theorem 3.4.11]{Lo}).

\begin{theorem}\label{B.d.thm}
In the assumptions and under the identifications of
Theorem~\ref{hkr}, the Connes-Tsygan differential
$$
B:HH_i(A) \to HH_{i+1}(A), \qquad i \geq 0
$$
becomes de Rham differential $d:\Omega^i_A \to
\Omega^{i+1}_A$.\endproof
\end{theorem}

\subsection{Cartier isomorphism.}\label{car.subs}

Now assume that $k$ is a perfect field of positive characteristic
$\cchar k = p$. Then in the assumptions of Theorem~\ref{hkr}, the
classic Cartier isomorphism gives an expression for the de Rham
cohomology groups of the variety $X = \Spec A$. Namely, denote by
$$
B\Omega_A^\hdot \subset Z\Omega_A^\hdot \subset \Omega^\hdot_A
$$
the image resp.\ kernel of the de Rham differential. Then one has a
functorial isomorphism
\begin{equation}\label{car.c}
C:Z\Omega^\hdot_A/B\Omega_A^\hdot \cong \Omega^\hdot_{A^{(1)}},
\end{equation}
where $A^{(1)}$ is the twist of $A$ by the Frobenius endomorphism $k
\to k$ (that is, $A$ with the $k$-vector space structure given by
$\lambda \cdot a = \lambda^pa$, $\lambda \in k$, $a \in A$). We will
need a non-commutative generalization of this given in \cite{cart}.

For any $k$-vector space $E$, denote
$$
C_{(1)}(E) = (E^\otimes p)_\sigma, \qquad C^{(1)}(E) = (E^{\otimes
  p})^\sigma,
$$
where $\sigma$ is the permutation of order $p$, and let $E^{(1)}$ be
the Frobenius twist of $E$. Then one has a natural map
\begin{equation}\label{C.E}
C:E^{(1)} \to C_{(1)}(E), \qquad e \mapsto e^{\otimes p},
\end{equation}
and one observes that the map is additive and $k$-linear. If $E$ is
finite-dimensional, then one can dualize the construction and obtain
a functorial $k$-linear map
\begin{equation}\label{R.E}
R:C^{(1)}(E) \to E^{(1)}.
\end{equation}
Taking filtered colimits, one extends the map $R$ to arbitrary
vector spaces. The cokernel $\Phi(E)$ of the map \eqref{C.E}
coincides with the kernel of the map \eqref{R.E}, and the
composition
$$
\begin{CD}
C_{(1)}(E) @>>> \Phi(E) @>>> C^{(1)}(E)
\end{CD}
$$
of the natural projection and the natural embedding is the trace map
$\tr_{\Z/p\Z}$ for the group $\Z/p\Z$ generated by $\sigma$.

If one now takes a $k$-algebra $A$ and considers the object $A^\hush
\in \Fun(\Lambda,k)$, then the maps $C$, $R$ taken together define
functorial maps
\begin{equation}\label{C.eq}
C:A^{(1)\hush} \to \pi_{p!}i_p^*A^\hush, \qquad
R:\pi_{p*}i_p^*A^\hush \to A^{(1)\hush}.
\end{equation}
The cokernel $\Phi A^\hush$ of the map $C$ coincides with the kernel
of the map $R$, and the composition map
$$
\begin{CD}
\pi_{p!}i_p^*A^\hush @>>> \Phi A^\hush @>>> \pi_{p*}i_p^* A^\hush
\end{CD}
$$
is the trace map $\tr_p$ pf \eqref{e.tr.pi}. One then considers the
complex $\K^p_\idot(A^\hush)$ of Definition~\ref{K.l.defn} and the
augmentation maps
$$
\begin{CD}
\pi_{p*}i_p^*A^\hush[1] @>{\wt{\kappa}_1}>>
\K^p(A^\hush) @>{\wt{\kappa}_0}>>
\pi_{p!}i_p^*A^\hush
\end{CD}
$$
induced by the exact sequence \eqref{K.l.4}, and one defines
subcomplexes $B\K^p_\idot(A^\hush)$ and $Z\K^p_\idot(A^\hush)$ in
$\K^p_\idot(A^\hush)$ by
\begin{equation}\label{B.Z.K}
B\K^p_\idot(A^\hush) = \wt{\kappa}_1(\Ker R), \qquad
Z\K^p_\idot(A^\hush) = \wt{\kappa}_0^{-1}(\Im C).
\end{equation}
Note that since $\wt{\kappa}_0 \circ \wt{\kappa}_1=0$, we have
$B\K^p_\idot(A^\hush) \subset Z\K^p_\idot(A^\hush) \subset
\K^p_\idot(A^\hush)$. In fact, $B\K^p_\idot(A^\hush) \cong \Phi
A^\hush[-1]$, and we have a short exact sequence
$$
\begin{CD}
0 @>>> Z\K^p_\idot(A^\hush) @>>> \K^p_\idot(A^\hush) @>>> \Phi
A^\hush @>>> 0.
\end{CD}
$$
The quotient $Z\K^p_\idot(A^\hush)/B\K^p_\idot(A^\hush)$ is a
complex of length $2$, and its homology objects in degrees $0$ and
$1$ are identified with $A^{(1)\hush}$ by the maps \eqref{C.eq}. In
this, it is similar to the complex $\K_\idot(A^{(1)\hush})$.

If the algebra $A = k[G]$ is monomial, then it is easy to identify
the two complexes. Namely, the diagonal maps $G^n \to G^{np}$, $n
\geq 1$ commute with the structure maps \eqref{a.hash.maps} and
define a map of functors $\pi_p^*G^\hush \to i_p^*G^\hush$. By
\eqref{A.G}, it induces maps $\pi_p^*A^{(1)\hush} \to i_p^*A^\hush$,
$\K_\idot(\pi_p^*A^{(1)\hush}) \to \K_\idot(i_p^*A^\hush)$, and by
adjunction, we obtain a map from $\K_\idot(A^{(1)\hush})$ to
$\K^p_\idot(A^\hush)$. By \cite[Lemma 5.2]{cart}, it lands in
$Z\K^p_\idot(A^\hush) \subset \K^p_\idot(A^\hush)$, and it is easy
to check that the induced map
\begin{equation}\label{Phi}
\K_\idot(A^{(1)\hush}) \to Z\K^p_\idot(A^\hush)/B\K^p_\idot(A^\hush)
\end{equation}
is a quasiisomorphism. For a general algebra $A$, there no obvious
map of complexes like this. Nevertheless, we have the following.

\begin{theorem}\label{car.thm}
Assume given a perfect field $k$ of characteristic $p=\cchar k > 0$
and an associative unital $k$-algebra $A$.
\begin{enumerate}
\item Denote by $\nu:\D(\Lambda,k) \to \D(\Lambda,W_2(k))$ the
  tautological functor that takes a $k$-vector spaces and treats it
  as a module over the second Witt vectors ring $W_2(k)$. Then there
  exists a canonical isomorphism
\begin{equation}\label{car.nc}
\nu(\K_\idot(A^{(1)\hush})) \cong
\nu(Z\K^p_\idot(A^\hush)/B\K^p_\idot(A^\hush))
\end{equation}
in the derived category $\D(\Lambda,W_2(k))$. If $A=k[G]$ is
monomial, then this isomorphism coincides with \eqref{Phi}.
\item Assume that $A$ is commutative and satisfies the assumptions
  of Theorem~\ref{hkr}. Then the isomorphisms \eqref{hkr.iso} and
  \eqref{K.l.cor} induce isomorphisms
$$
HC_i(B\K^p_\idot(A^\hush)) \cong B\Omega^i_A, \qquad
HC_i(Z\K^p_\idot(A^\hush)) \cong Z\Omega^i_A,
$$
and with these identifications, the map induced by the map
\eqref{car.nc} is inverse to the Cartier isomorphism
\eqref{car.c}.
\end{enumerate}
\end{theorem}

\proof[Start of the proof.] If the characteristic $p$ is odd, then
the isomorphism \eqref{car.nc} is \cite[Lemma 4.4]{cart}, it
coincides with \eqref{Phi} by \cite[Lemma 5.2]{cart}, and there is
no need to apply the functor $\nu$. For $p=2$, one has to apply
$\nu$ and modify the argument as in \cite[Section 4]{dege}. This
gives \eqref{car.nc} but does not prove its compatibility with
\eqref{Phi}. Moreover, we will need a more explicit form of
\eqref{car.nc}, and for this, we need some notions to be introduced
later. Thus we postpone the end of the proof of
Theorem~\ref{car.thm}~\thetag{i} to the end of
Subsection~\ref{basic.subs}. As for \thetag{ii}, it is
\cite[Proposition 5.1]{cart}. While formally, that paper requires
$p$ to be odd throughout, this is not really used in Section~5 that
contains Propo\-si\-tion~5.1 --- all one needs is \cite[Lemma
  5.2]{cart} that reduces to the compatibility between
\eqref{car.nc} and \eqref{Phi}.
\endproof

Finally, in Section~\ref{ill.sec}, we will need one simple result
that morally belongs to \cite{cart} but is not contained explicitly
in that paper; we prove it here.

\begin{lemma}\label{free.le}
For any associative unital $k$-algebra $A$ and integer $n \geq 1$,
the object $i_{p^n}^*\Phi A^\hush$ is $\pi_{p^n}$-free.
\end{lemma}

\proof{} By definition, it suffices to prove that for any $k$-vector
space $E$, the $\Z/p^n\Z$-action on $\Phi(E^{\otimes p^n})$ induced by
the $\Z/p^{n+1}\Z$-action on $E^{\otimes p^{n+1}}$ is free. Choose a
basis in $E$, so that $E = k[S]$ for some set $S$, and decompose
$$
E^{\otimes p^{n+1}} = k[S^{p^{n+1}}] = E^{\otimes p^n} \oplus E',
$$
where $E^{\otimes p^n} = k[S^{p^n}]$ is embedded by the diagonal
embedding $S^{p^n} \subset S^{p^{n-1}}$, and $E'=k[S']$ is spanned
by the complement $S'=S^{p^{n+1}} \setminus S^{p^n}$. Then we have
$$
\Phi(E^{\otimes p^n}) \cong E'_{\sigma^{p^n}},
$$
where $\sigma$ is the generator of the group $\Z/p^{n+1}\Z$, and
this group acts on $S'$ without fixed points.
\endproof

\subsection{Trace functors.}\label{trace.subs}

It is useful to categorify the construction of the cyclic object
$A^\hush$ in the following way. Assume given a unital monoidal
category $\C$, and define a category $\C^\hush$ as follows.
\begin{itemize}
\item Objects are pairs $\langle [n],\{c_y\}\rangle$ of an object
  $[n] \in \Lambda$ and a collection of objects $c_y$ in $\C$
  numbered by elements $y \in Y([n])$.
\item Morphisms from $\langle [n'],\{c'_y\}\rangle$ to $\langle
  [n],\{c_y\}\rangle$ are given by a collection of a morphism
  $f:[n'] \to [n]$ and morphisms
$$
f_y:\bigotimes_{y' \in f^{-1}(y)}c'_{y'} \to c_y, \qquad y \in
Y([n]),
$$
where the product is taken in the canonical order on $f^{-1}(y)$.
\end{itemize}
The forgetful functor $\rho:\C^\hush \to \Lambda$ is a
cofibration. Its fiber $\C^\hush_{[n]}$ over $[n] \in \Lambda$ is
naturally identified with the product $\C^{Y([n])}$ of copies of
$\C$ numbered by elements $y \in Y([n])$, and the transition
functors are induced by the tensor product in $\C$ via an obvious
generalization of \eqref{a.hash.maps}. A morphism $\langle f,\{f_y\}
\rangle$ in $\C^\hush$ is cocartesian with respect to $\rho$ if and
only if all its components $f_y$ are invertible maps.

\begin{defn}\label{trace.def}
A {\em trace functor} from $\C$ to some category $\E$ is a
functor $F^\hush:C^\hush \to \E$ that sends all maps cocartesian
with respect to $\rho:\C^\hush \to \Lambda$ to invertible maps in
$\E$.
\end{defn}

Definition~\ref{trace.def} is a version of \cite[Definition
  2.1]{trace}, with the equivalence proved in \cite[Lemma
  2.5]{trace}. More precisely, every trace functor $F^\hush:\C^\hush
\to \E$ induces a functor $F:\C \to \E$ by restriction to $\C =
\C^\hush_{[1]} \subset \C^\hush$. Then \cite[Lemma 2.5]{trace} shows
that extending an arbitrary $F:\C \to \E$ to a trace functor is
equivalent to giving functorial isomorphisms
\begin{equation}\label{tau.mn}
\tau_{M,N}:F(M \otimes N) \cong F(N \otimes M), \qquad M,N \in \C
\end{equation}
that satisfy some compatibility constraints listed in
\cite[Definition 2.1]{trace}.

Any associative unital algebra object $A$ in a unital monoidal
category $\C$ canonically defines a section
\begin{equation}\label{alpha}
\alpha:\Lambda \to \C^\hush
\end{equation}
of the projection $\rho$. On objects, $\alpha$ sends $[n] \in
\Lambda$ to $\langle [n],\{A\} \rangle$ (that is, we take $n$ copies
of $A$). On morphisms, $\alpha$ is essentially given by
\eqref{a.hash.maps} (for a more invariant categorical definition, see
\cite[Subsection 1.4]{trace}). Given a trace functor
$F^\hush:\C^\hush \to \E$ to some category $\E$, we can compose it
with $\alpha$ and obtain a canonical functor
\begin{equation}\label{FA.eq}
FA^\hush = F^\hush \circ \alpha:\Lambda \to \E.
\end{equation}
For any $[n] \in \Lambda$, we can choose a morphism $f:[n] \to [1]$,
and its cocartesian lifting proveds an isomorphism
\begin{equation}\label{FA.n}
FA^\hush([n]) \cong F(A^{\otimes n}).
\end{equation}
Note, however, that the left-hand side carries a canonical action of
$\Z/n\Z=\Aut([n])$, and the right-hand side does not: the symmetry
is broken by the choice of the morphism $f$. To recover the
$\Z/n\Z$-action on $F(A^{\otimes n})$, one has to use the maps
\eqref{tau.mn} (specifically, $\tau_{A,A^{\otimes n-1}}$ generates
the group).

If the monoidal category $\C$ is symmetric (for example, $\C$ is the
category of flat modules over a commutative ring $k$), then the
tautological functor $I:\C \to \C$ carries natural isomorphisms
\eqref{tau.mn} provided by the symmetry map; these satify the
necessary constraints and promote $I$ to a trace functor
$I^\hush:\C^\hush \to \C$. Then $A^\hush = IA^\hush$ is exactly the
canonical cyclic object of Subsection~\ref{alg.subs}.

We will also need a version of the isomorphism \eqref{kun.loc} for
the objects $FA^\hush$. To obtain it, recall that a {\em
  pseudotensor structure} on a functor $F:\C \to \E$ between unital
monoidal categories $\C$, $\E$ with unit objects $1_\C$, $1_\E$ and
product functors $m_\C:\C^2 \to \C$, $m_\E:\E^2 \to \E$ is given by
functorial maps
\begin{equation}\label{eps.mu}
\eps:1_\E \to F(1_\C), \ \mu_{M,N}:F(M \otimes N) \to F(M) \otimes
F(N), \quad M,N \in \C,
\end{equation}
such that $\mu_{1,M}=\eps \otimes \id$, $\mu_{M,1}=\id \otimes \eps$
(unitality), and $(\id \otimes \mu_{N,L}) \circ \mu_{M,N \otimes L}
= (\mu_{M,N} \otimes \id) \circ \mu_{M \otimes N,L}$ for any $M,N,L
\in \C$ (associativity). If $\C$ and $\E$ are symmetric monoidal
categories, then a pseudotensor structure is {\em symmetric} if
$\mu$ commutes with commutativity isomorphisms. Now assume given a
unital monoidal category $\E$ and a symmetric unital monoidal
category $\C$. Then $\C^2 = \C \times \C$ is also a unital monoidal
category, and since $\C$ is symmetric, the tensor product functor
$m_{\C}:\C^2 \to \C$ is a tensor functor. Therefore it extends to a
functor
$$
m^\hush:\C^{2\hush} \cong \C^\hush \times_\Lambda \C^\hush \to \C
$$
that commutes with projections $\rho$. Moreover, this relative
tensor product $m^\hush$ is symmetric in the obvious sense, and the
unit object $1_\C \in \C$ is an algebra object, so that we obtain a
section $1^\hush:\Lambda \to \C^\hush$ of the projection
$\rho$. This section is obviously cocartesian.

\begin{defn}\label{tr.ten.def}
A {\em pseudotensor structure} on a trace functor $F^\hush:\C^\hush
\to \E$ is given by functorial maps
\begin{equation}\label{eps.mu.hash}
\eps:1_{\E} \to F^\hush(1^\hush), \qquad \mu:m_{\E} \circ F^{2\hush}
\to F^\hush \circ m^\hush
\end{equation}
satisfying the unitality and associativity conditions. If $\E$ is
symmetric, then a pseudotensor structure is {\em symmetric} if $\mu$
commutes with the commutativity isomorphisms.
\end{defn}

For any pseudotensor structure on a trace functor $F^\hush$, the
functor $F:\C \to \E$ induced by $F^\hush$ inherits a pseudotensor
structure in the usual sense. Conversely, a pseudotensor structure
on $F$ extends to $F^\hush$ if it is compatible with the maps
\eqref{tau.mn}.

Now for any two associative unital algebra objects $A$, $B$ in a
symmetric monoidal category $\C$, and for any trace functor
$F^\hush:\C^\hush \to \E$ equpped with a pseudotensor structure in
the sense of Definition~\ref{tr.ten.def}, we obtain a canonical map
\begin{equation}\label{ab.loc}
\mu:FA^\hush \otimes FB^\hush \to F(A \otimes B)^\hush,
\end{equation}
a version of the isomorphism \eqref{kun.loc}. If the pseudoensor
structure on $F^\hush$ is tensor --- that is, if the maps
\eqref{eps.mu.hash} are isomorphisms --- then \eqref{ab.loc} is also an
isomorphism. This is the case, for instance, when $F^\hush$ is the
tautological functor $I^\hush$. If $\E$ and the pseudotensor
structure on $F^\hush$ are symmetric, then the map \eqref{ab.loc}
commutes with the symmetry maps. If, moreover, $A$ is commutative,
then $FA^\hush$ becomes a commutative associative unital algebra
object in the category of functors from $\Lambda$ to $\E$.

\subsection{Polynomial Witt vectors.}\label{witt.subs}

Assume now given a perfect field $k$ of positive characteristic
$p$. For any integer $m \geq 1$, let $W_m(k)$ be the ring of
$m$-truncated $p$-typical Witt vectors of $k$, and let $W(k)$ be
their inverse limit. Then \cite[Section 2]{witt} constructs a series
of polynomial functors
$$
W_m:k\amod \to W_m(k)\amod, \qquad m \geq 1
$$
and surjective restriction maps $R:W_{m+1} \to W_m$, $m \geq 1$,
called {\em polynomial Witt vectors} functors. It is further proved
in \cite[Section 3]{witt} that one has injective co-restriction maps
$C:W_m \to W_{m+1}$, $m \geq 1$ such that $C \circ R = p \id$, and
non-additive but functorial Teichm\"uller maps
\begin{equation}\label{T.W}
T:E \to W_m(E), \qquad E \in k\amod
\end{equation}
such that $R^{m-1} \circ T = \id$. Moreover, $W_m$ together with the
maps $R$, $C$ and $T$ canonically extend to trace functors
$W_m^\hush:k\amod^\hush \to W_m(k)\amod$.

To describe the construction, recall from \cite[Subsection
  4.1]{witt} that for any unital mono\-i\-dal category $\C$ and
integer $l \geq 1$, the functors $i_l$, $\pi_l$ of \eqref{i.pi.l}
give rise to a natural commutative diagram
\begin{equation}\label{big.cart}
\begin{CD}
\C^\hush @<{i^{\C}_l}<< \C^\hush_l @>{\pi^{\C}_l}>> \C^\hush\\
@V{\rho}VV @VV{\rho}V @VV{\rho}V\\
\Lambda @<{i_l}<< \Lambda_l @>{\pi_l}>> \C^\hush,
\end{CD}
\end{equation}
where the square on the right-hands side is Cartesian, and the
functor $i^{\C}_l$ sends $\langle [n],c_\idot \rangle$ to $\langle
i_l([n]),c^l_\idot \rangle$, with the collection $c^l_\idot$ given
by $c^l_y = c_{\eta_l(y)}$, $y \in Y(i_l([n]))$, where $\eta_l$ is
the map \eqref{q.l}.

We now fix an integer $m \geq 1$, and let $\C_m$ be the category of
flat finitely generated $W_m(k)$-modules, with the natural quotient
functor
$$
q:\C_m \to \C, \qquad E \mapsto E/p
$$
to the category $\C=\C_1$ of finite-dimensional $k$-vector
spaces. Moreover, let $l=p^m$, and simplify notation by writing
$i^{(m)} = i^{\C_m}_{p^m}$, $\pi^{(m)}=\pi^{\C_m}_{p^m}$, with
$\vpi^{(m)}_*$ being the corresponding functor \eqref{v.ga}. Then
$\C_m$ is a symmetric unital tensor category, and the tautological
embedding $\C_m \subset W_m(k)\amod$ defines a canonical trace
functor $I^\hush \in \Fun(\C^\hush_m,W_m(k))$. By base change, the
functor
$$
Q^\hush_m = \vpi^{(m)}_*i^{(m)*}I^\hush \in \Fun(\C_m^\hush,W_m(k))
$$
is also a trace functor. Here is, then, the main result of
\cite{witt}.

\begin{theorem}[{{\cite[Proposition~4.3~\thetag{i}]{witt}}}]\label{witt.thm}
There exists a unique trace functor $W_m^\hush \in \Fun(\C,W_m(k))$
such that $Q_m^\hush \cong q^*W_m^\hush$.\endproof
\end{theorem}

This completely defines $W_m^\hush$ on finite-dimensional $k$-vector
spaces; one then extends it to all vector spaces by requiring that
the extension commutes with filtered colimits. In addition,
\cite[Proposition~4.3~\thetag{ii}]{witt} constructs the restriction
maps
\begin{equation}\label{R.eq}
R:W^\hush_{m+1} \to W^\hush_m
\end{equation}
and the Teichm\"uller maps $T:W^\hush_1 \to W^\hush_m$, and
\cite[Lemma~3.1]{witt} provides co-restriction maps $C:W^\hush_m \to
W^\hush_{m+1}$.  Passing to the limit with respect to the
restriction maps, one obtaines a trace functor $W^\hush$ from
$k$-vector spaces to modules over the Witt vectors ring $W(k)$.

By definition, the functors $W_m^\hush$ carry a decreasing {\em
  standard filtration} $F^iW_m^\hush$ and an increaing {\em
  co-standard filtration} $F_iW_m^\hush$ given by
$$
F^iW_m^\hush = \Ker R^{m-i},F_iW_m^\hush = \Im C^{m-i} \subset
W_m^\hush,
$$
where $R$ and $C$ are the restriction and the co-restriction maps.
In the limit $W^\hush$, only the standard filtration survives, and
we have $W_m^\hush = W^\hush/F^mW^\hush$.  It has been also proved
in \cite[Subsection 3.1]{witt} that in fact the values of the limit
functor $W^\hush$ are torsion-free, and for any $k$-vector space
$E$, we have
$$
W(E)/p \cong \lim_{\overset{R}{\gets}}C^{(i)}(E),
$$
where $C^{(i)}(E) \subset E^{\otimes p^i}$ is the subspace of
$\Z/p^i\Z$-invariant vectors.

If $m=1$, then $W_1^\hush$ is simply the tautological trace functor
$I^\hush$. In particular, it carries a pseudotensor structure. This
has been extended to all $m \geq 1$ in
\cite[Proposition~4.7]{witt}. Namely, it has been proved that the
functors $W^\hush_m$, $m \geq 1$ carry natural pseudotensor
structures in the sense of Definition~\ref{tr.ten.def}, and the
restriction maps $R:W^\hush_{m+1} \to W^\hush_m$ together with the
Teichm\"uller maps \eqref{T.W} are compatible with the pseudotensor
structures. Passing to the limit, one obtains a pseudotensor
structure on the trace functor $W^\hush$.

Moreover, consider the diagram \eqref{big.cart} for $l=p$, and
simplify notaton by writing $i=i^\C_p$, $\pi=\pi^\C_p$. Then we have
the following.

\begin{prop}\label{FV.prop}
For any $m \geq 1$, there exists functorial additive maps
$$
V:W^\hush_m \circ i \to W^\hush_{m+1} \circ \pi, \qquad
F:W^\hush_{m+1} \circ \pi \to W^\hush_m \circ i
$$
such that $RV=VR$, $FR=RF$, $VF=p\id$, and $FV$ coincides with the
trace map $\tr^\dg_\pi$ of \eqref{tr.dg}.
\end{prop}

We can again pass to the limit with respect to the maps $R$ and
obtain maps $V$, $F$ for the limit trace functor $W^\hush$. They
define an $FV$-structure on $W^\hush$ in the sense of
\cite[Definition~4.5]{witt}, and by \cite[Proposition~4.7]{witt},
this structure is compatible with the pseudotensor structure in the
sense of \cite[Definition~4.6]{witt}.

\section{Definitions and properties.}\label{defs.sec}

\subsection{General case.}\label{basic.subs}

As in Subsection~\ref{witt.subs}, fix a perfect base field $k$ of
positive characteristic $p$. Consider the trace functors
$W_m^\hush$, $m \geq 1$ provided by Theorem~\ref{witt.thm}, the
restriction maps $R$ of \eqref{R.eq}, the Teichm\"uller maps $T$ of
\eqref{T.W}, and the co-restriction maps $C$ such that
$RC=CR=p\id$. Then for any associative unital $k$-algebra $A$ and
integer $m \geq 1$, we have a natural object
\begin{equation}\label{W.m.A}
W_mA^\hush \in \Fun(\Lambda,W_m(k))
\end{equation}
of \eqref{FA.eq}, and we have natural maps
\begin{equation}\label{RC.WA}
R:W_{m+1}A^\hush \to W_mA^\hush, \qquad C:W_mA^\hush \to
W_{m+1}A^\hush
\end{equation}
and the Teichm\"uller map
\begin{equation}\label{T.A}
A^\hush \to W_mA^\hush
\end{equation}
of simplicial sets. We also have the limit trace functor $W^\hush$
and the limit object
\begin{equation}\label{W.A}
WA^\hush \cong \lim_{\overset{R}{\gets}}W_mA^\hush
\end{equation}
in the category $\Fun(\Lambda,W(k))$. For any $m$, $W_mA^\hush$
carries the standard and costandard filrations $F^\idot W_mA^\hush$,
$F_\idot W_mA^\hush$ induced by the corresponding filtrations on the
trace functor $W_m^\hush$. The standard filtration survives in the
limit, and we have $W_mA^\hush = WA^\hush/F^mWA^\hush$ for any $m
\geq 1$. Moreover, $WA^\hush$ is torsion-free, and we have
$$
WA^\hush/p = \lim_{\overset{R}{\gets}}\pi_{p^n*}i_{p^n}^*A^\hush,
$$
where $R$ is the map \eqref{C.eq}.  If $m=1$, then
$W_1^\hush=I^\hush$ is the tautological trace functor, so that
$W_1A^\hush \cong A^\hush$.

\begin{defn}
The {\em Hochschild-Witt homology groups} of the algebra $A$ are
given by
$$
W_mHH_\idot(A) = HH_\idot(W_mA^\hush), \quad m \geq 1, \qquad
WHH_\idot(A) = HH_\idot(WA^\hush),
$$
where $W_mA^\hush$, $WA^\hush$ are the natural cyclic objects
\eqref{W.m.A}, \eqref{W.A}, and $HH_\idot(-)$ are the Hochschild
homology groups of Definition~\ref{hh.defn}. The {\em
  Hochschild-Witt complex} $WCH_\idot(A)$ and its truncated versions
$W_mCH_\idot(A)$ are given by
$$
W_mCH_\idot(A) = CH_\idot(W_mA^\hush), \quad m \geq 1, \qquad
WCH_\idot(A) = CH_\idot(WA^\hush),
$$
where $CH_\idot(-)$ is the Hochschild complex of
Definition~\ref{hh.defn}.
\end{defn}

Since the trace functors $W^\hush_m$, $m \geq 1$ are equipped with
pseudotensor structures, \eqref{ab.loc} gives a natural map
\begin{equation}\label{W.AB}
W_mA^\hush \otimes W_mB^\hush \to W_m(A \otimes B)^\hush, \qquad m
\geq 1
\end{equation}
for any two associative unital $k$-algebras $A$, $B$, and these maps
are compatible with the restriction maps $R$. They are also
compatible with the Teichm\"uller maps \eqref{T.A}, in the sense
that by \cite[(3.16)]{witt}, we have a commutative diagram
\begin{equation}\label{T.mult}
\begin{CD}
A^\hush \times B^\hush @>>> (A \otimes B)^\hush\\
@V{T \times T}VV @VV{T}V\\
W_mA^\hush \times W_mB^\hush @>>> W_m(A \otimes B)^\hush
\end{CD}
\end{equation}
of cyclic sets. In the limit, we obtain a map
$$
WA^\hush \otimes WB^\hush \to W(A \otimes B)^\hush.
$$
Neither this map nor the maps \eqref{W.AB} are isomorphisms (unless
$m=1$). However, coupled with the K\"unneth map \eqref{kunn}, they
do provide product maps
\begin{equation}\label{WHH.prod}
\begin{aligned}
W_mHH_\idot(A) \otimes W_mHH_\idot(A) &\to W_mHH_\idot(A \otimes B),
\quad m \geq 1,\\
WHH_\idot(A) \otimes WHH_\idot(B) &\to WHH_\idot(A \otimes B)
\end{aligned}
\end{equation}
on the level of Hochschild-Witt homology; these maps are associative
and unital, and compatible with the restriction maps. For a
commutative algebra $A$, $WHH_\idot(A)$ and $W_mHH_\idot(A)$, $m
\geq 1$ thus become unital graded-commutative algebra, and the
Connes-Tsygan differential $B$ is its derivation.

On the other hand, note that by definition, the functor
$\alpha:\Lambda \to k\amod^\hush$ of \eqref{alpha} corresponding to
the algebra $A$ induces a functor $\alpha_p:\Lambda_p \to
k\amod^\hush_p$ such that
$$
\alpha \circ \pi_p \cong \pi_p^{k\amod} \circ \alpha_p, \qquad
\alpha \circ i_p \cong i_p^{k\amod} \circ \alpha_p.
$$
Then Proposition~\ref{FV.prop} provides natural functorial maps
\begin{equation}\label{FV.WA}
V:i_p^*W_mA^\hush \to \pi_p^*W_{m+1}A^\hush, \qquad
F:\pi_p^*W_{m+1}A^\hush \to i_p^*W_mA^\hush
\end{equation}
for any $m \geq 1$, and in the limit, we obtain a natural
$FV$-structure on the cyclic abelian group $WA^\hush$. On the level
of homology, the maps \eqref{FV.WA} induce natural group maps
\begin{equation}\label{FV.WHH}
\begin{aligned}
&V:W_mHH_\idot(A) \to W_{m+1}HH_\idot(A),\\
&F:W_{m+1}HH_\idot(A) \to W_mHH_\idot(A)
\end{aligned}
\end{equation}
for any $m \geq 1$, and the corresponding maps in the limit. By
definition, we have $VR=RV$, $FR=RF$, $VC=CV$, $FC=CF$, where $R$
and $C$ are the maps \eqref{RC.WA}.

\begin{lemma}\label{FV.prod.le}
For any integer $m \geq 1$ and two unital associative $k$-algebras
$A$, $B$, we have
$$
\begin{aligned}
F(a \cdot b) &= F(a) \cdot F(b), \quad a \in W_{m+1}HH_\idot(A), b
\in W_{m+1}HH_\idot(B),\\
V(a \cdot F(b)) &= V(a) \cdot b, \qquad\;\; a \in W_mHH_\idot(A), b \in
W_{m+1}HH_\idot(B),\\
V(F(b) \cdot a) &= a \cdot V(b), \qquad\;\;  a \in W_{m+1}HH_\idot(A), b \in
W_mHH_\idot(B),
\end{aligned}
$$
where the product is taken with respect to the maps
\eqref{FV.WHH}. If $A$ is commutative, and we take $B=A$, then all
three equalities hold with respect to the product in the algebras
$W_\idot HH_\idot(A)$.
\end{lemma}

\proof{} Immediately follows from \cite[Definition~4.6]{witt}.
\endproof

Moreover, we can iterate the maps \eqref{FV.WA} as in
\eqref{FV.ite}. In particular, for any integers $m, n \geq 1$, we
obtain natural maps of cyclic abelian groups
\begin{equation}\label{bar.FV.WA}
\overline{V}^n:\pi_{p^n!}i_{p^n}^*W_mA^\hush \to W_{m+n}A^\hush,
\quad
\overline{F}^n:W_{m+n}A^\hush \to \pi_{p^n*}i_{p^n}^*W_mA^\hush,
\end{equation}
a version of the maps \eqref{bar.FV.ite}.

\begin{prop}\label{whh.prop}
For any integers $m, n \geq 1$ and any associative unital
$k$-algebra $A$, the maps \eqref{FV.ite} and \eqref{RC.WA} fit into
short exact sequences of cyclic abelian groups
$$
\begin{CD}
0 @>>> \pi_{p^n!}i_{p^n}^*W_mA^\hush @>{\overline{V}^n}>>
W_{m+n}A^\hush @>{R^m}>> W_nA^\hush @>>> 0,\\
0 @>>> W_nA^\hush @>{C^m}>> W_{m+n}A^\hush @>{\overline{F}^n}>>
\pi_{p^n*}i_{p^n}^*W_mA^\hush @>>> 0.
\end{CD}
$$
Moreover, for any $m \geq 1$, we have commutative diagrams
$$
\begin{CD}
\pi_{p^{m-1}!}i_{p^{m-1}}^*A^{\hush} @>{C}>>
\pi_{p^m!}i_{p^m}^*A^\hush\\
@V{\overline{V}^{m-1}}VV @VV{\overline{V}^m}V\\
W_mA^\hush @>{C}>> W_{m+1}A^\hush
\end{CD}
\qquad
\begin{CD}
\pi_{p^m*}i_{p^m}^*A^\hush @>{R}>>
\pi_{p^{m-1}*}i_{p^{m-1}}^*A^{\hush}\\
@A{\overline{F}^m}AA @AA{\overline{F}^{m-1}}A\\
W_{m+1}A^\hush @>{R}>> W_mA^\hush
\end{CD}
$$
where $C$ and $R$ in the top row are obtained by applying
$\pi_{p^{m-1}!}i_{p^{m-1}}^*$ resp.\ $\pi_{p^{m-1}*}i_{p^{m-1}}^*$
to the maps $C$ resp.\ $R$ of \eqref{C.eq}. 
\end{prop}

\proof{} It suffices to prove everything after evaluation at an
arbitrary object $[n] \in \Lambda$. Choosing a map $f:[n] \to [1]$
provides an identification $W_mA^\hush([n]) \cong W_m(A^{\otimes
  n})$, and by \cite[Proposition~4.7~\thetag{iii}]{witt}, our
sequences then become the exact sequences of \cite[Lemma~3.7]{witt}
for $E=A^{\otimes m}$. The commutative diagrams are those of
\cite[Lemma~3.4]{witt}.
\endproof

\begin{corr}\label{gr.corr}
The associated graded quotients $\gr^\hdot_F$ and $\gr_\idot^F$ of
the objects $WA^\hush$, $W_mA^\hush$ with respect to the standard
and co-standard filtrations are given by
$$
\gr^i_FW_mA^\hush \cong \gr^i_FWA^\hush
\cong \pi_{p^{i-1}!}i_{p^{i-1}}^*A^\hush, \quad
\gr_i^FW_mA^\hush \cong \pi_{p^{i-1}*}i_{p^{i-1}}^*A^\hush.
$$
\end{corr}

\proof{} Clear. \endproof

\begin{lemma}\label{fdv.whh.le}
For any associative unital $k$-algebra $A$ and integer $m \geq 1$,
we have
$$
FV=VF=p\id:W_mHH_\idot(A) \to W_mHH_\idot(A),
$$
and
$$
FBV=B:W_mHH_\idot(A) \to W_mHH_{\idot+1}(A),
$$
where $F$ and $V$ are the maps \eqref{FV.WHH}, and $B$ is the
Connes-Tsygan differential. The same equalities hold in the limit
groups $WHH_\idot(A)$.
\end{lemma}

\proof{} Everything except for the equality $VF=p\id$ is
Lemma~\ref{fdv.loc}. The equality $VF=p\id$ follows from
\cite[Lemma~3.11]{witt} by the same argument as in the proof of
Proposition~\ref{whh.prop}.
\endproof

\proof[End of the proof of Theorem~\ref{car.thm}~\thetag{i}.] For
any complex $E_\idot$ in an abelian category $\E$, and any integers
$n \leq m$, denote by $\tau_{[n,m]}E_\idot$ the truncation of
$E_\idot$ with respect to the canonical filtration --- that is, let
$\ol{E}_\idot = \tau_{[n,m]}E_\idot$ be the subquotient of $E_\idot$
given by
$$
\ol{E}_n = \Ker d_n, \qquad \ol{E}_m = \Coker d_{m+1}, \qquad \ol{E}_i=E_i,
\quad n < i < m,
$$
and $\ol{E}_i=0$ otherwise, where $d_i:E_i \to E_{i-1}$ is the
differential. Recall that the functors $\tau_{[n,m]}$ descend to the
derived category $\D(\E)$.

For any ring $O$ and any object $E \in \Fun(\Lambda_p,O)$, extend
the resolution \eqref{reso} to the right to obtain a $2$-periodic
complex $P_\idot(E)$ in $\Fun(\Lambda_p,O)$ with terms
$$
K_{2i}(E) = \K_0(E), \quad K_{2i+1}(E) = \K_1(E), \qquad i \in \Z,
$$
and denote by $\pi_{p\flat}E \in \D(\Lambda,O)$ the object
represented by the complex $\pi_{p*}P_\idot(E)$. Then
\cite[(3.10)]{dege} provides a functorial map
\begin{equation}\label{a.eq}
\K_\idot(\Z)[1] \otimes (\tau_{[-1,0]}\pi_{p\flat}E)[1] \to
\tau_{[0,1]}\pi_{p\flat}E.
\end{equation}
Moreover, for any object $E' \in \D(\Lambda,O)$ equipped with a map
$b:E' \to (\tau_{[-1,0]}\pi_{p\flat}E)[1]$, we can compose $b$ with the
map \eqref{a.eq} and obtain a map
\begin{equation}\label{bb.eq}
\K_\idot(\Z) \otimes E' \to \tau_{[-1,0]}\pi_{p\flat}E.
\end{equation}
If the composition map
\begin{equation}\label{comp.bb}
\begin{CD}
E' @>{b}>> (\tau_{[-1,0]}\pi_{p\flat}E)[1] @>>>
(\tau_{[-1,-1]}\pi_{p\flat}E)[1]
\end{CD}
\end{equation}
is an isomorphism, then \eqref{bb.eq} is an isomorphism in
$\D(\Lambda,O)$.

Now take $O=W_2(k)$ and $E=\nu(i_p^*A^\hush)$. Then by definition,
we have an isomorphism of complexes
$$
\nu(Z\K^p_\idot(A^\hush)/B\K^p_\idot(A^\hush)) \cong
\tau_{[0,1]}\pi_{p*}P_\idot(E).
$$
On the other hand, the shifted truncation
$(\tau_{[-1,0]}\pi_{p*}P_\idot(E))[1]$ is the complex $\wt{E}_\idot$
with terms $\wt{E}_1=\pi_{p!}E$, $\wt{E}_0=\pi_{p*}E$, and with the
differential given by the trace map $\tr_p$.

To obtain the isomorphism \eqref{car.nc}, one takes as $E'$ the
complex of length $2$ given by $E'_1=\wt{E}_1$, $E'_0=W_2 A^\hush$,
and with the differential given by the map $V$ of
\eqref{FV.WA}. Then we have a map $b:E'_\idot \to \wt{E}_\idot$
equal to $\id$ in degree $1$ and $F$ in degree $0$, and the
composition
$$
\begin{CD}
E'_\idot @>{b}>> \wt{E}_\idot @>>> \tau_{[0,0]}\pi_{p\flat}E \cong
\nu(A^{(1)\hush})
\end{CD}
$$
is a quasiisomorphism by Proposition~\ref{whh.prop}. The
corresponding map \eqref{bb.eq} is the isomorphism \eqref{car.nc}.

To see that this isomorphism is compatible with \eqref{Phi}, note
that for any set $S$, sending an element $s \in S$ to its
Teichm\"uller represetative $T(s)$ gives a natural map $W_2(k)[S]
\to W_2(k[S])$, and if we compose it with the projection
$F:W_2(k[S]) \to C^{(1)}(k[S])$, then by \cite[(3.18)]{witt}, the
resulting map is induced by the diagonal embedding $S \subset
S^p$. If $A=k[G]$, we can collect these maps for all powers $G^n$,
$n \geq 1$, and obtain a natural map
\begin{equation}\label{e.e.eq}
E''_\idot \to E'_\idot,
\end{equation}
where $E''_\idot$ is given by $E''_0 = W_2(k)[G^\hush]$,
$E''_1=k[G^\hush]$, with the differential induced by the embedding
$V:k \to W_2(k)$. The map \eqref{e.e.eq} is a quasiisomorphism, and
the map \eqref{bb.eq} for the object $E''_\idot$ exactly coincides
with the map \eqref{Phi}.
\endproof

\subsection{Free algebras.}\label{free.subs}

Now we want to give some idea as to how Hochschild-Witt homology
groups look like. In order to do this, we compute them for free
algebras over $k$. So, let $M$ be a $k$-vector space, and let
$$
A = T^\hdot(M) = \bigoplus_{n \geq 0}M^{\otimes n}
$$
be the free associative unital $k$-algebra generated by $M$. By
definition, the algebra $A$ is graded by non-negative integers, and
its Hochschild homology $HH_\idot(A)$ inherits the grading. Recall
that for any integer $n \geq 1$, the component $HH_\idot(A)_n
\subset HH_\idot(A)$ of degree $n$ is given by
\begin{equation}\label{hh.ten}
HH_i(A)_n \cong
\begin{cases}
(M^{\otimes n})_\sigma, &\quad i = 0,\\
(M^{\otimes n})^\sigma, &\quad i = 1,\\
0, &\quad i \geq 2,
\end{cases}
\end{equation}
where $\sigma:M^{\otimes n} \to M^{\otimes n}$ is the permutation of
order $n$ that generates the action of the cyclic group
$\Z/n\Z$. The component $HH_i(A)_0$ is $k$ for $i=0$ and $0$
otherwise. The Connes-Tsygan differential $B$ vanishes on
$HH_0(A)_0$, and for $n \geq 1$, we have
$$
B = \tr_{\Z/n\Z}:(M^{\otimes n})_\sigma \to (M^{\otimes n})^\sigma,
$$
where the $\Z/n\Z$-action is generated by $\sigma$

To see the isomorphism \eqref{hh.ten} explicitly, one can use the
functors $e^n_!$ of \eqref{e.n}. For any $n \geq 1$, the map $M \to
A$ induces a $\Z/n\Z$-equivariant map $M^{\otimes n} \to A^{\otimes
  n}_n = e^{n*}A^\hush_n$, where $A^\hush_n \subset A^\hush$ is the
component of degree $n$, and by adjunction, we obtain a natural map
$$
e^n_!M^{\otimes n} \to A^\hush_n,
$$
An easy computation shows that this map is in fact an isomorphism,
so that we have
$$
HH_\idot(e^n_!M^{\otimes n}) \cong HH_\idot(A)_n.
$$
Then Lemma~\ref{e.n.le} immediately implies \eqref{hh.ten}.

In order to generalize \eqref{hh.ten} to Hochschild-Witt homology,
recall from \cite[Subsection~3.4]{witt} that a $\Z$-grading on a
$k$-vector space $E$ tautologically induces a $\Z$-grading on
$W_m(E)$ for any $m \geq 1$, but these gradings are not compatible
with the restriction maps $R$. In order to make them compatible, one
has to rescale them, so that the natural grading on $W_m(E)$ and on
the limit group $W(E)$ is indexed by elements $a \in \Z[1/p]$ in the
localization at $p$ of the ring $\Z$. Thus $W_mHH_\idot(A)$ and
$W_mCH_\idot(A)$ are also graded by $\Z[1/p]$. Note that the
Frobenius and Verschibung maps $F$, $V$ of \eqref{FV.WHH} do not
preserve the rescaled degree: $F$ has degree $p$, while $V$ has
degree $1/p$.

For any $n \in \Z$, $m \geq 1$, we can equip $W_m(M^{\otimes n})$
with a $\Z/n\Z$-action by \eqref{FA.n}, and we then have a
$\Z/n\Z$-equivariant map $W_m(M^{\otimes n}) \to
e^{n*}W_mA^\hush_n$, where as before, $W_mA^\hush_n$ stands for the
component of degree $n$. By adjunction, we get a natural map
\begin{equation}\label{e.n.WA}
e^n_!W_m(M^{\otimes n}) \to W_mA^\hush_n.
\end{equation}
If $m \geq 2$, this map is no longer an isomorphism. However, we
still have the following result.

\begin{lemma}\label{free.pos.le}
For any $n \geq 1$, the map
$$
\K_\idot(W_m(M^{\otimes n})) \to W_mCH_\idot(A)_n
$$
induced by \eqref{e.n.WA} and Lemma~\ref{e.n.le} is a
quasiisomorphism.
\end{lemma}

\proof{} It suffices to prove that the map \eqref{e.n.WA} induces a
quasiisomorphism on Hochschild homology of the associated graded
quotients $\gr^i_F$ with respect to the standard filtration. Since
the functor $e^n_!$ is exact, we have $\gr^i_Fe^n_!W_m(M^{\otimes
  n}) \cong e^n_!\gr^i_FW_m(M^{\otimes n})$, while the quotient
$\gr^i_FW_mCH_\idot(A)$ is described by
Corollary~\ref{gr.corr}. With these identifications, the claim
immediately follows from Proposition~\ref{Y.prop}.
\endproof

Assume now given $a \in \Z[1/p]$ of the form $a = np^{-i}$, $n > 0$
prime to $p$, $i > 0$. We then have the following result.

\begin{lemma}\label{free.neg.le}
For any $m \geq 1$, the iterated Verschiebung and Frobenius maps
\begin{equation}\label{v.i.f.i}
\begin{aligned}
V^i:&W_mHH_0(A)_n \to W_{m+i}HH_0(A)_a,\\
F^i:&W_{m+i}HH_1(A)_a \to W_mHH_1(A)_n
\end{aligned}
\end{equation}
are isomorphisms, and $W_mHH_i(A)_a=0$ for $i \neq 0,1$.
\end{lemma}

\proof{} Denote $q=p^i$. Corollary~\ref{gr.corr} and
Proposition~\ref{whh.prop} immediately show that the map
$\overline{V}^i:\pi_{q!}i_q^*W_mA^\hush_m \to
W_{m+i}A^\hush_a$ becomes an isomorphism after taking the associated
graded quotients with respect to the standard filtration. Therefore
it is itself an isomorphism, and together with
Lemma~\ref{free.pos.le}, we obtain a quasiisomorphism
\begin{equation}\label{w.m.ch}
W_{m+i}CH_\idot(A)_a \cong
CH_\idot(\pi_{q!}i_q^*e^n_!W_m(M^{\otimes n})).
\end{equation}
For any $\Z[\Z/n\Z]$-module $E$, denote
$K^{n,q}_\idot(E)=CH_\idot(i_q^*e^n_!E)$. Then explicitly, we have
$$
i_q^*e^n_!E([m]) \cong E[Y^{n,q}([m])]_{\Z/n\Z}, \qquad [m]
\in \Lambda_q,
$$
where $Y^{n,q}:\Lambda_q \to \Sets$ is the functor \eqref{Y.nm}. By
Lemma~\ref{fp.le}, the $\Z/q\Z$-action on $Y^{n,q}([m])$ is free for
any $[m] \in \Lambda_q$, so that $i_q^*e^n_!E$ is
$\pi_q$-free. Moreover, $Y^{n,q} \circ j^o_q:\Delta^o \to \Sets$ is
a simplicial set with a free $\Z/q\Z$-action and only finitely many
non-degenerate simplices. Therefore the Tate cohomology
$\vH^\hdot(\Z/q\Z,K^{n,q}_\idot(E))$ vanishes, and we have natural
quasiisomorphisms
\begin{equation}\label{L.R.M}
\begin{aligned}
C_\idot(\Z/q\Z,K^{n,q}_\idot(E)) &\cong K_\idot^{n,q}(E)_{\Z/q\Z}
\overset{\tr}{\cong}\\
&\overset{\tr}{\cong} K_\idot^{n,q}(E)^{\Z/q\Z} \cong
H^\hdot(\Z/q\Z,K_\idot^{n,q}(E)),
\end{aligned}
\end{equation}
where $C_\idot(\Z/q\Z,-)$, $C^\hdot(\Z/q\Z,-)$ are the complexes
computing homology and cohomology of the group $\Z/q\Z$. If we take
$E = W_m(M^{\otimes n})$, then $K_\idot^{n,q}(E)$ is quasiisomorphic
to $W_mCH_\idot(A)_n$ by Lemma~\ref{free.pos.le} and
Lemma~\ref{la.ho}, \eqref{w.m.ch} identifies the quasiisomorphic
complexes \eqref{L.R.M} with $W_{m+1}CH_\idot(A)_a$, and the maps
$V^i$ resp.\ $F^i$ of \eqref{v.i.f.i} are induced by the natural
adjunction maps
$$
K_\idot^{n,q}(E) \to C_\idot(\Z/q\Z,K_\idot^{n,q}(E)), \qquad
C^\hdot(\Z/q\Z,K_\idot^{n,q}(E)) \to K_\idot^{n,q}(E).
$$
But the complex $K_\idot^{n,q}(E)$ only has non-trivial homology groups in
degrees $0$ and $1$, and the $\Z/q\Z$-action on these groups is
trivial by Corollary~\ref{i.l.triv}. This proves the claim.
\endproof

We can now combine Lemma~\ref{free.pos.le} and
Lemma~\ref{free.neg.le} to obtain a complete description of the
Hochschild-Witt homology groups of the algebra $A = T^\hdot(M)$. For
simplicity, we only treat the limit case $WHH_\idot(A)$ and leave it
to the reader to figure out the truncated versions. We note that any
positive element $a \in \Z[1/p]$, $a > 0$ can be uniquely
represented as $a = np^i$, where $n > 0$ is prime to $p$, and we
denote $i=\|a\|$.

\begin{theorem}\label{free.thm}
The Hochschild-Witt homology $WHH_i(A)$ vanishes if $i \neq 0,1$,
and $WHH_i(A)_a=0$ for negative $a \in \Z[1/p]$. For $a=0$, we have
$WHH_0(A)_0=k$ and $WHH_1(A)_0=0$. For any $a > 0$ with $\|a\| \geq 0$,
we have natural identifications
\begin{equation}\label{free.1}
WHH_0(A)_a \cong W(M^{\otimes a})_{\Z/a\Z}, \qquad
WHH_1(A)_a \cong W(M^{\otimes a})^{\Z/a\Z},
\end{equation}
where $\Z/a\Z$ acts on $W(M^{\otimes a})$ via the isomorphism
\eqref{FA.n}. The Frobenius map $F:WHH_1(A)_a \to WHH_1(A)_{pa}$ is
an isomorphism for any $a > 0$. The Verschiebung map
$V:WHH_0(A)_{pa} \to WHH_0(A)_a$ is an isomorphism if $\|a\| < 0$.
If $\|a\| \geq 0$, so that $a = np^{\|a\|}$ is an integer, the
iterated Verschiebung map
\begin{equation}\label{free.2}
V^{\|a\|}:WHH_0(A)_a \to WHH_0(A)_n
\end{equation}
is injective and identifies $WHH_0(A)_a$ with the image of the
iterated Verschiebung map
\begin{equation}\label{free.3}
\overline{V}^{\|a\|}:W(M^{\otimes a})_{\Z/a\Z} \to W(M^{\otimes
  n})_{\Z/n\Z} = WHH_0(A)_n.
\end{equation}
\end{theorem}

\proof{} We note that the information contained in
Theorem~\ref{free.thm} together with the identities of
Lemma~\ref{fdv.whh.le} allow one to reconstruct immediately the maps
$F$ and $V$ in all the other degrees, and also compute the
Connes-Tsygan differential $B$. In particular, in the situation when
isomorphisms \eqref{free.1} hold, $B$ must coincide with the trace
map $\tr_{\Z/a\Z}$. This also immediately follows from
Lemma~\ref{free.pos.le} and Lemma~\ref{e.n.le}, as indeed do the
isomorphisms \eqref{free.1}. The fact that $V$ and $F$ are
isomorphisms when $\|a\| < 0$ follows from Lemma~\ref{free.neg.le}
by taking the limit with respect to $m$. Finally, the fact that the
maps \eqref{free.2} and \eqref{free.3} agree immediately follows
from the construction of the quasiisomorphisms of
Lemma~\ref{free.pos.le}, and since $W(E)$ has no $p$-torsion for any
$E$, Lemma~\ref{fdv.whh.le} together with \cite[Lemma 3.11]{witt}
then implies that the Frobenius map $F^{\|a\|}$ coincides with the
corresponding Frobenius map $\ol{F}^{\|a\|}$. The latter is an
isomorphism by \cite[Corollary 3.9]{witt}.
\endproof

\section{Comparison I.}\label{hess.sec}

We now want to compare the Hochschild-Witt homology $WHH_0(A)$ of
degree $0$ of an associative unital algebra $A$ over a perfect field
$k$ of characteristic $p=\cchar k > 0$ to the group $W(A)$ of
non-commutative Witt vectors of $A$ constructed in \cite{hewi}. We
start by reviewing Hesselholt's construction.

\subsection{Recollection.}\label{hess.subs}

Hesselholt starts with an arbitrary associative ring $A$, possibly
non-unital, and defines the {\em ghost map} $w:A^{\N} \to
(A/[A,A])^{\N}$ by
\begin{equation}\label{w.i}
w(a_0 \times a_1 \times \dots ) = w_0 \times w_1 \times \dots,
\qquad w_i = \sum_{j=0}^i p^ja_j^{p^{i-j}},
\end{equation}
where $A/[A,A]=HH_0(A)$ is the quotient of $A$ by the abelian
subgroup spanned by all commutators $aa'-a'a$, $a,a' \in A$. One
observes that $w$ is well-defined as a map of sets. One also
observes that for any $n \geq 0$, $w_n$ only depends on $a_i$ with
$0 \leq i \leq n$, so that we also have truncated ghost maps
$w:A^{n+1} \to (A/[A,A])^{n+1}$, where the terms in the product
correspond to the integers $0,\dots,n \in \N$. If we denote by
$R:A^{n+1} \to A^n$ the restriction map obtained by projecting to
the first $n$ coordinates, then $w \circ R = R \circ w$.

With these data, Hesselholt constructs a series of abelian groups
$W_n(A)$, $n \geq 1$, and surjective set-theoretic maps $q:A^n \to
W_n(A)$ such that
\begin{enumerate}
\renewcommand{\labelenumi}{(\Roman{enumi})}
\item the ghost map $w:A^n \to (A/[A,A])^n$
  factors as
$$
\begin{CD}
A^n @>{q}>> W_n(A) @>{\overline{w}}>> (A/[A,A])^n,
\end{CD}
$$
where $\overline{w}$ is a map of abelian groups,
\item $W_n(A)$ and $\overline{w}$ are functorial in $A$,
\item if $A/[A,A]$ has no $p$-torsion, then $\overline{w}$ is
  injective, and
\item we have functorial surjective restriction map $R:W_{n+1}(A)
  \to W_n(A)$ such that $R \circ q = q \circ R$ and $\overline{w}
  \circ R = R \circ \overline{w}$.
\end{enumerate}
The group $W(A)$ is then the inverse limit of $W_n(A)$ with respect
to $R$.

The construction of the groups $W_n(A)$ given in \cite{hewi} was not
quite correct, so let us present a corrected construction given in
\cite{heerr}. The construction is inductive. For $n=1$, we just take
$W_1(A) = A/[A,A]$, $\overline{w}=\id$, and we observe that all the
properties are tautologically satisfied. We then fix $n \geq 1$ and
assume given groups $W_i(A)$ as above for all $i \leq n$. The
construction proceeds in two steps. In the first step, consider the
set
$$
\wt{W}_{n+1}(A) = A \times W_n(A)
$$
and the map
$$
\wt{q} = \id \times q:A^{n+1} = A \times A^n \to \wt{W}_{n+1}(A),
$$
where the first component in $A \times A^n$ corresponds to the
$0$-th component in $A^{n+1}$. Denote by $\wt{R}:\wt{W}_{n+1}(A) \to
A$ the projection onto this first component, so that we have $\wt{R}
\circ \wt{q} = R^n$. Note that by \eqref{w.i} and the condition
\thetag{I}, the ghost map $w:A^{n+1} \to (A/[A,A])^{n+1}$ factors as
$w = \wt{w} \circ \wt{q}$, with $\wt{w}:\wt{W}_{n+1}(A) \to
(A/[A,A])^{n+1}$ given by
\begin{equation}\label{wt.w}
\wt{w}(a \times b) = w(a \times 0)+(0 \times p\overline{w}(b)),
\qquad a \in A, b \in W_n(A).
\end{equation}

\begin{lemma}\label{c.i.le}
There exist universal polynomials $c_i(s_0,s_1)$, $i \geq 1$ of
degrees $p^i$ in two non-commuting variables $s_0$, $s_1$ such that
for any $n \geq 1$ and elements $a_0,a_1 \in A$ in an associative
ring $A$, we have
\begin{equation}\label{c.i.f}
(a_0+a_1)^{p^n}-a_0^{p^n}-a_1^{p^n}=\sum_{i=1}^np^ic_i(a_0,a_1)^{p^{n-i}}
\mod [A,A].
\end{equation}
\end{lemma}

\proof{} It obviously suffices to consider the universal case
$A=T^\hdot(M)$, the tensor algebra generated by the free abelian
group $M = \Z[\{s_0,s_1\}]$. Then \cite[Lemma~3.15~\thetag{i}]{witt}
provides elements $c_i \in M^{\otimes p^i}$ such that for any $n
\geq 1$, we have
$$
(s_0+s_1)^{\otimes p^n} = s_0^{\otimes p^n} + s_1^{\otimes p^n} +
  \sum_{i=1}^n\sum_{j=0}^{p^i-1}\sigma^j(c_i^{\otimes p^{n-i}}) \in
  M^{\otimes p^n},
$$
where for any $l \geq 0$, $\sigma:M^{\otimes l} \to M^{\otimes l}$ is the
cyclic permutation. To deduce the claim, it remains to recall that
for any $m \in M^{\otimes l} \subset A$, the difference $m -
\sigma(m)$ lies in the subspace $[A,A] \subset A$.

Alternatively, one can use the polynomials $\delta_i$ of
\cite{hewi}; these work by \cite[Proposition~1.2.3]{hewi}.
\endproof

\begin{lemma}\label{W.le}
There exists a unique functorial abelian group structure on
$\wt{W}_{n+1}(A)$ such that
\begin{enumerate}
\item $\wt{R}:\wt{W}_{n+1}(A) \to A$ is a group map, and
\item $\wt{w}:\wt{W}_{n+1}(A) \to (A/[A,A])^{n+1}$ is also a group
  map.
\end{enumerate}
Moreover, in this group structure, we have
\begin{equation}\label{coc.W}
(a \times b) + (a' \times b') = (a+a') \times (b+b'-q(c_1(a,a')
\times \dots \times c_n(a,a')))
\end{equation}
for any $a,a' \in A$, $b,b' \in W_n(A)$, where $c_\idot(-,-)$ are
the polynomials provided by Lemma~\ref{c.i.le}.
\end{lemma}

\proof{} In principle, existence follows from \cite{heerr}. Namely,
Hesselholt constructs a binary operation on $A^{n+1}$ given by
explicit universal non-commu\-ta\-tive polynomials
$s_i(a_0,\dots,a_n,a_0',\dots,a'_n)$, $0 \leq i \leq n$, and then
shows that it descends to an abelian group structure on
$\wt{W}_{n+1}A$. He then checks that $w:A^{n+1} \to (A/[A,A])^{n+1}$
factors as $w = \wt{w} \circ \wt{q}$ for some group map
$\wt{w}:\wt{W}_{n+1}(A) \to (A/[A,A])^{n+1}$. Since $\wt{q}$ is
surjective, the factorization is unique, so that Hesselholt's
$\wt{w}$ must coincide with the map \eqref{wt.w}. This gives
\thetag{ii}. Moreover, it is immediately clear from his construction
that $s_i$ does not depend on $n$ and only depends on
$a_0,\dots,a_i,a_0',\dots,a_i'$, and in particular,
$s_0=a_0+a_0'$. This is \thetag{i}.

Alternatively, it is easy to see existence directly. Namely, one can
simply define the group structure on $\wt{W}_{n+1}(A)$ by
\eqref{coc.W}. Then \thetag{i} is given, and \thetag{ii} immediately
follows from \eqref{w.i}, \eqref{wt.w} and \eqref{c.i.f}. The only
thing to check is that we indeed have an abelian group structure ---
that is, that the map
\begin{equation}\label{coc.loc}
A \times A \to W_n(A), \qquad a \times a' \mapsto q(c_1(a,a') \times
\dots \times c_n(a,a')))
\end{equation}
is a symmetric $2$-cocycle of the group $A$ with values in
$W_n(A)$. For every ring $A$, we can find a ring $A'$ and a
surjective map $A' \to A$ such that $A'/[A',A']$ has no $p$-torsion,
and since the map \eqref{coc.loc} is functorial, it suffices to
check the symmetic cocycle condition after replacing $A$ with
$A'$. But then the ghost map $\overline{w}:W_n(A) \to (A/[A,A])^n$
is injective by \thetag{III}, and by \eqref{wt.w}, $\wt{w}$ is also
injective on $a \times W_n(A) \subset \wt{W}_{n+1}(A)$ for any $a
\in A$. Therefore it suffices to check the symmetric cocycle
condition after applying the map $\wt{w}$, and then it immediately
follows from \thetag{ii}.

Uniqueness follows from \eqref{coc.W}, and to prove \eqref{coc.W},
note that by \thetag{i}, the group law in $\wt{W}_{n+1}(A)$ is in
any case given by
$$
(a \times b) + (a' \times b') = (a+a') \times (b+b'-F(a,a',b,b')),
$$
where $F$ is a certain functorial map
$$
F:A \times A \times W_n(A) \times W_n(A) \to W_n(A).
$$
We need to prove that $F(a,a',b,b')=q(c_1(a,a') \times \dots \times
c_n(a,a'))$. As before, by functoriality, it suffices to prove this
claim for rings $A$ such that $A/[A,A]$ has no $p$-torsion, and in
this case, $\wt{w}$ is injective on $W_n(A) \times a \subset
\wt{W}_{n+1}(A)$ for any $a \in A$. Thus it suffices to prove the
claim after applying the map $\wt{w}$, and then it immediately
follows from Lemma~\ref{c.i.le}.
\endproof

\begin{remark}
While we need Lemma~\ref{W.le} for comparison purposes, the reader
may notice that it can also be used as an alternative for
Hesselholt's construction of the group structure on
$\wt{W}_{n+1}(A)$ and of the map $\wt{w}$. It is essentially the same
argument, but it needs a smaller number of formulas.
\end{remark}

In the second step of the construction, Hesselholt considers the
free abelian group $\Z[A \times A]$ spanned by the set $A \times A$,
and the map
\begin{equation}\label{wt.d}
\wt{d}:\Z[A \times A] \to \wt{W}_{n+1}(A), \qquad \wt{d}(a \times
a') = (aa'-a'a) \times 0.
\end{equation}
One then sets $W_{n+1}(A) = \Coker \wt{d}$. By \eqref{wt.w}, the
composition $\wt{w} \circ \wt{d}$ vanishes, so that $\wt{w}$ factors
through a functorial map
$$
\overline{w}:A^{n+1} \to W_n(A).
$$
This gives \thetag{I}, and \thetag{II}, \thetag{IV} are also obvious
from the construction. To finish the inductive step, one has to
check \thetag{III}. There seems to be no direct argument for this
statement, so Hesselholt proves it by topological methods (and also
proves a comparison theorem between $W_n(A)$ and a certain
functorial group $TR^n_0(A;p)$ that comes from algebraic topology).

One additional thing that comes out of Hesselholt's construction is
a functorial map $V:W_n(A) \to W_{n+1}(A)$ induced by the embedding
$0 \times W_n(A) \subset A \times W_n(A) = \wt{W}_{n+1}(A)$. By
definition, it fits into a functorial exact sequence
\begin{equation}\label{hess.ex}
\begin{CD}
W_n(A) @>{V}>> W_{n+1}(A) @>{R^n}>> A @>>> 0.
\end{CD}
\end{equation}
In general, the sequence need not be exact on the left (and an
example when it is not is contained in \cite{heerr}). Hesselholt
further observes that the map \eqref{wt.d} factors through the
quotient $A \otimes A$ of the group $\Z[A \times A]$, and that it
then induces a functorial map $\partial:HH_1(A) \to W_n(A)$ that
extends the exact sequence \eqref{hess.ex} one step to the left.

\subsection{Comparison.}

Assume now given an associative unital algebra $A$ over a perfect
field $k$ of characteristic $p=\cchar k > 0$, and consider the
Hochschild-Witt homology groups $W_nHH_0(A)$, $n \geq 1$. By
definition, we have
$$
W_nHH_0(A) = H_0(\Delta^o,j^{o*}W_nA^\hush) =
\colim_{\Delta^o}j^{o*}W_nA^\hush,
$$
so that by adjunction, we have a natural augmentation map
$\tau:j^{o*}W_nA^\hush \to W_nHH_0(A)$, where $W_nHH_0(A)$ is
understood as a constant simplicial group. We also have the
Teichm\"uller map $T:A^\hush \to W_nA^\hush$ of
\eqref{T.A}. Composing the two maps and evaluating at $[1] \in
\Delta^o$, we obtain a natural map of sets
\begin{equation}\label{bar.q}
\overline{q} = (\tau \circ T)([1]):A \to W_nHH_0(A).
\end{equation}
In addition, we have the restriction maps $R:W_nHH_0(A) \to
W_{n-1}HH_0(A)$ and the Verschiebung maps $V:W_{n-1}HH_0(A) \to
W_nHH_0(A)$ of \eqref{FV.WHH}. Combining the maps $V$ with
\eqref{bar.q}, we obtain a functorial map of sets
\begin{equation}\label{q.eq}
q:A^n \to W_nHH_0(A), \qquad q(a_1 \times \dots \times a_n) =
\sum_{i=1}^n V^{i-1}(\overline{q}(a_i)).
\end{equation}
To avoid confusion with the polynomial Witt vectors $W_n$ of
Subsection~\ref{witt.subs}, let us from now denote Witt vectors of
Subsection~\ref{hess.subs} by $W_n^H(A)$ (where $H$ stands for
``Hesselholt''). The comparion theorem that we want to prove is the
following one.

\begin{theorem}\label{hess.thm}
There exists functorial isomorphisms
$$
\iota:W_nHH_0(A) \cong W^H_n(A),  \qquad n \geq 1
$$
such that $R \circ \iota = \iota \circ R$, $V \circ \iota = \iota
\circ V$, and $\iota \circ q = q$.
\end{theorem}

In order to prove this, we need to find Hochschild-Witt counterparts
of all the steps in the inductive construction presented in
Subsection~\ref{hess.subs}. For $n=1$, Theorem~\ref{hess.thm} is
clear: we have $W_1HH_0(A) = HH_0(A) = A/[A,A]$, and this coincides
with $W_1^H(A)$ by definition. Take an integer $n \geq 2$. Denote
$W_n^1A^\hush = \pi_{p!}i_p^*W_{n-1}A^\hush$ and
$$ 
W_n^1CH_\idot(A) = CH_\idot(W^1_nA^\hush),\qquad 
W_n^1HH_\idot(A) = HH_\idot(W^1_nA^\hush).
$$
Then by Proposition~\ref{whh.prop}, we have a natural short exact
sequence
\begin{equation}\label{w.n1.ex}
\begin{CD}
0 @>>> W_n^1A^\hush @>{\overline{V}}>> W_nA^\hush
@>{R^{n-1}}>> A^\hush @>>> 0
\end{CD}
\end{equation}
in $\Fun(\Lambda,W_n(k))$ that induces an exact sequence of
Hochschild homology complexes. By Lemma~\ref{la.ho}, we also have
$$
\begin{aligned}
W_n^1HH_0(A) &= \colim_{\Delta^o}j^{o*}\pi_{p!}i_p^*W_{n-1}A^\hush
\cong \colim_{\Delta^o \times \ppt_p}j^{o*}_pi_p^*W_{n-1}A^\hush
\cong\\
&\cong \colim_{\ppt_p}W_{n-1}HH_0(A) = W_{n-1}HH_0(A)_{\Z/p\Z},
\end{aligned}
$$
and since the $\Z/p\Z$-action on $W_{n-1}HH_0(A)$ is trivial by
Corollary~\ref{i.l.triv}, this group is identified with
$W_{n-1}HH_0(A)$. Let $\tau:j^{o*}W_n^1A^\hush \to W_{n-1}HH_0(A)$
be the augmentation map, and let $\bW_nA^\hush =
j^{o*}W_nA^\hush/\overline{V}(\Ker \tau)$. Then \eqref{w.n1.ex}
induces a short exact sequence
\begin{equation}\label{bar.w.ex}
\begin{CD}
0 @>>> W_{n-1}HH_0(A) @>>> \bW_nA^\hush @>>> j^{o*}A^\hush @>>> 0
\end{CD}
\end{equation}
of simplicial $W_n(k)$-modules, where $W_{n-1}HH_0(A)$ is understood
as a constant simplicial module. Thus if we denote by
$$
\bW_nCH_\idot(A) = C_\idot(\bW_nA^\hush)
$$
the standard complex of the simplicial $W_n(k)$-module
$\bW_nA^\hush$, then in degree $0$, we have a natural short exact
sequence
\begin{equation}\label{wt.whh.ex}
\begin{CD}
0 @>>> W_{n-1}HH_0(A) @>{V}>> \bW_nCH_0(A) @>{R^{n-1}}>> A
@>>> 0,
\end{CD}
\end{equation}
where $A=CH_0(A)$ is the degree-$0$ term of the Hochschild complex
$CH_\idot(A)$. The Teichm\"uller map $T$ provides a set-theoretic
splitting $T:j^{o*}A^\hush \to \bW_nA^\hush$ of the sequence
\eqref{bar.w.ex} and an isomorphism of sets
\begin{equation}\label{wt.whh.spl}
\overline{W}_nCH_0(A) \cong A \times W_{n-1}HH_0(A).
\end{equation}

\begin{lemma}\label{whh.wt.le}
In terms of the isomorphism \eqref{wt.whh.spl}, the group law in the
abelian group $\bW_nCH_0(A)$ is given by
\begin{equation}\label{coc.whh}
(a \times b) + (a' \times b') = (a+a') \times (b+b'-q(c_1(a,a')
\times \dots \times c_n(a,a'))),
\end{equation}
where $q$ is the map \eqref{q.eq}, and $c_i$ are the polynomials of
Lemma~\ref{c.i.le}.
\end{lemma}

\proof{} By the uniqueness clause of Lemma~\ref{W.le}, it suffices
to prove the claim for one particular choice of the polynomials
$c_i$; let us choose those provided by
\cite[Lemma~3.15~\thetag{i}]{witt}.  For any $k$-vector space $E$,
\cite[Lemma~3.7]{witt} provides a short exact sequence
$$
\begin{CD}
0 @>>> \bW^1_n(E) @>{\overline{V}}>> W_n(E)
@>{R^{n-1}}>> E @>>> 0,
\end{CD}
$$
where $\bW^1_n(E) = W_{n-1}(E^{\otimes p})_{\Z/p\Z}$, and $\Z/p\Z$
acts on $W_{n-1}(E^{\otimes p})$ via the isomorphism
\eqref{FA.n}. The Teichm\"uller map $T:E \to W_n(E)$ splits this
sequence set-theoretically, so that the group law in $W_n(E) =
E \times \bW^1_n(E)$ is given by
\begin{equation}\label{coc.E}
(a \times b) + (a' \times b') = (a+a') \times (b+b'-\wt{c}(a,a'))
\end{equation}
for any $a,a' \in E$, $b,b' \in \bW^1_n(E)$, where $\wt{c}:E \times
E \to \bW^1_n(E)$ is a certain $2$-cocycle of the group $E$ with
coefficients in $\bW^1_n(E)$. The cocycle $\wt{c}$ has been computed
in \cite[Lemma~3.15~\thetag{ii}]{witt} --- one has
\begin{equation}\label{wt.c}
\wt{c}(e,e') = \sum_{m=1}^nV^m(T(c_m(e,e'))), \qquad e,e' \in E,
\end{equation}
where $T:E^{\otimes p^m} \to W_{n-m}(E^{\otimes p^m})$ is the
Teichm\"uller splitting, the polynomials $c_m$ are evaluated in the
tensor algebra $T^\hdot E$, so that $c_m(e,e')$ lies in $E^{\otimes
  p^m}$, and $V^m:W_{n-m}(E^{\otimes p^m}) \to W_n(E)$ is the
iteration of the Verschiebung map. Now observe that for any $m \geq
1$, we have a commutative diagram of simplicial sets
$$
\begin{CD}
\bi_{p^m}^*j^{o*}A^\hush @>{T}>> \bi_{p^m}^*j^{o*}W_{n-m}A^\hush\\
@V{c_{p^m}}VV @VV{c_{p^m}}V\\
j^{o*}A^\hush @>{T}>> j^{o*}W_{n-m}A^\hush
\end{CD}
$$
and a commutative diagram of simplicial $W_n(k)$-modules
$$
\begin{CD}
\bi_{p^m}^*j^{o*}W_{n-m}A^\hush @>{V^m}>> j^{o*}W_nA^\hush
@>{\tau}>> W_nHH_0(A)\\ 
@VV{c_{p^m}}V @. @|\\
j^{o*}W_{n-m}A^\hush @>{\tau}>> W_{n-m}HH_0(A) @>{V^m}>> W_nHH_0(A),
\end{CD}
$$
where $T$ are the Teichm\"uller maps, $c_{p^m}$ are the egdewise
subdivision maps \eqref{h.l}, and $\tau$ are the augmentation
maps. It remains to notice that after evaluation at $[1] \in
\Lambda$, the map
$$
c_{p^m}:A^{\otimes p^m} = \bi_{p^m}^*j^{o*}A^\hush([1]) \to A =
j^{o*}A^\hush([a])
$$
sends $a_1 \otimes \dots \otimes a_{p^m}$ to the product $a_1\dots
a_{p^m} \in A$, and compare \eqref{coc.whh}, \eqref{bar.q} and
\eqref{q.eq} with \eqref{coc.E} and \eqref{wt.c}.
\endproof

\proof[Proof of Theorem~\ref{hess.thm}.] Fix an integer $n \geq 1$,
and assume by induction that the isomorphisms $\iota$ are already
constructed for $W_iHH_0(A)$ with $i \leq n$. Then comparing
Lemma~\ref{whh.wt.le} and Lemma~\ref{W.le}, we see that the
isomorphism $\iota:W_nHH_0(A) \cong W_n^H(A)$ canonically extends to
a functorial isomorphism
$$
\wt{\iota}:\bW_{n+1}CH_0(A) \cong \wt{W}_{n+1}(A)
$$
that automatically commutes with the maps $V$ and $R$. To show that
$\wt{\iota}$ then descends to an isomorphism $\iota$ between
$W_{n+1}HH_0(A)$ and $W_{n+1}^H(A)$, note that by definition, the
homology of the complex $\bW_{n+1}CH_\idot(A)$ in degree $0$ is
isomorphic to $W_{n+1}HH_0(A)$. The differential $d:\bW_{n+1}CH_1(A)
\to \bW_{n+1}CH_0(A)$ is given by $d = \partial_0 - \partial_1$,
where for $i=0,1$,
$$
\partial_i:\bW_{n+1}CH_1(A) = \bW_{n+1}A^\hush([2]) \to
\bW_{n+1}CH_0(A) = \bW_{n+1}A^\hush([1]),
$$
are the face maps corresponding to the two maps $[1] \to [2]$ in
$\Delta$. Then the Teichm\"uller map provides an isomorphism
$$
T:A^{\otimes 2} = A^\hush([2]) \cong \bW_{n+1}A^\hush([2]),
$$
and since $T$ is a map of simplicial sets, we have $T \circ
\partial_i = \partial_i \circ T$, $i=0,1$. The differential $d$
is therefore given by
$$
\begin{aligned}
d(a \otimes a') &= \partial_0(T(a \otimes a')) - \partial_1(T(a
\otimes a'))=\\
&= T(\partial_0(a \otimes a')) - T(\partial_1(a \otimes a')) =\\
&= T(aa') - T(a'a),
\end{aligned}
$$
for any $a \otimes a' \in A^{\otimes 2} \cong \bW_{n+1}CH_1(A)$, and
in terms of the identification \eqref{wt.whh.spl}, this coincides on
the nose with the map \eqref{wt.d}.
\endproof

\begin{remark}
Hesselholt's differential $\partial:HH_1(A) \to W_n^H(A)$ is also
visible in our approach --- this is simply the connecting
differential for the long exact sequence of Hochschild homology
corresponding to \eqref{bar.w.ex} (or \eqref{w.n1.ex}).
\end{remark}

\section{Comparison II.}\label{ill.sec}

To finish the paper, we now show that in the setup of the
Hochschild-Kostant-Rosenberg Theorem, our Hoshchild-Witt complex
reduces to de Rham-\-Witt complex of Deligne and Illusie \cite{ill}.

\subsection{Recollection.} We start by recalling the relevant
material from \cite{ill}. We fix a finite field $k$ of
characteristic $\cchar k = p$, and a commutative associative
$k$-algebra $A$. Recall that classically, one associates to $A$ the
commutative associative ring $W(A)$ of {\em $p$-typical Witt
  vectors} of $A$ (see e.g.\ \cite[Section 0.1]{ill} or any of the
other standard references). The construction is functorial, and
since $A$ is a $k$-algebra, $W(A)$ is a $W(k)$-algebra. As a set,
$W(A)$ is identified with the product $A^{\N}$ of a countable number
of copies of $A$, so that its elements are infinite sequences
$(x_0,x_1,\dots )$ of elements of $A$. The operations are given by
certain universal polynomials in $x_i$. One defines $n$-truncated
Witt vectors $W_n(A)$ by forgetting $x_i$, $i \geq n$; the
operations on $W(A)$ are compatible with the truncations, so that
$W_n(A)$ are also rings. We have natural additive multiplicative
restriction maps $R:W_n \to W_{n-1}(A)$. Moreover, one defines the
Verschiebung map $V:W_{n-1}(A) \to W_n(A)$ by $V((x_0,\dots,x_n)) =
(0,x_0,\dots,x_n)$, and the Frobenius map $F:W_{n+1}(A) \to W_n(A)$
by $F((x_0,\dots,x_n)) = (x_0^p,\dots,x_{n-1}^p)$. One has $FV = VF
= p$, $FR=RF$, $RV=VR$. For any $n,m \geq 1$, we have a natural
short exact sequence
$$
\begin{CD}
0 @>>> W_n(A) @>{V^m}>> W_{n+m}(A) @>{R^n}>> W_m(A) @>>>
0.
\end{CD}
$$
We have $W_1(A) \cong A$ canonically, and we have a compatible
system of natural multiplicative splittings $A \to W_n(A)$ of the
natural projections $R^n:W_n(A) \to A$; for any $x \in A$, we denote
by $\ux \in W_n(A)$ its image under the splitting, known as the
Teichm\"uller representative of $x$. This splitting is characterized
by the fact that $F\ux = \underline{x^p} = \ux^p$.

From now on, to avoid confusion with our polynomial Witt vectors, we
will denote Witt vectors ring $W_n(A)$, $W(A)$ by $W_n^{cl}(A)$,
$W^{cl}(A)$ (where $cl$ stands for ``classical''). Hesselholt's
non-commutative construction of \cite{hewi} is a generalization of
the classical construction, so that for a commutative associative
ring $A$, one has $W^H(A) \cong W^{cl}(A)$ and $W^H_n(A) \cong
W_n^{cl}(A)$, $n \geq 1$. These isomorphisms are compatible with the
restriction maps $R$ and the Verschiebung maps $V$.

\medskip

Now assume that as in Theorem~\ref{hkr}, $A$ is finitely
generated over $k$, and $X = \Spec A$ is smooth over $k$.

\begin{defn}[{{\cite[Definition I.1.1]{ill}}}]\label{V.pro}
A {\em $V$-de Rham procomplex} $M^\hdot_\idot$ on $X=\Spec A$ is a
collection $M^\hdot_n$, $n \geq 1$ of commutative associative DG
algebras and maps
$$
R:M_{n+1} \longrightarrow M_n, \qquad V:M_n \longrightarrow M_{n+1}
$$
between them such that $RV=VR$, $R$ is surjective, multiplicative
and commutes with the differential, and
\begin{enumerate}
\item the rings $M^0_n$ are idenitified with the rings $W_n^{cl}(A)$
  in such a way that $R$ and $V$ are the restriction and the
  Verschiebung maps,
\item we have
$$
V(xdy) = V(x) d V(y)
$$
for any $n$, $i$, $j$, $x \in M^i_n$, $y \in M^j_n$, and
\item we have
$$
V(y) d \ux = V(\ux^{p-1}y)dV(\ux)
$$
for any $n$, $x,y \in A$.
\end{enumerate}
\end{defn}

\begin{defn}\label{FV.pro}
An {\em $FV$-de Rham procomplex} $M^\hdot_\idot$ on $X$ is a
collection $M^\hdot_n$, $n \geq 1$ of commutative associative DG
algebras and maps
$$
R,F:M_{n+1} \longrightarrow M_n, \qquad V:M_n \longrightarrow M_{n+1}
$$
between them such that $RV=VR$, $RF=FR$, $FV=VF$, $R$ is
multiplicative and commutes with the differential, and
\begin{enumerate}
\item the rings $M^0_n$ are idenitified with the rings $W_n^{cl}(A)$
  in such a way that $R$, $F$ and $V$ are the restriction, the
  Frobenius and the Verschiebung maps, and
\item we have
$$
FV=VF=p, \qquad FdV=d,
$$
\item we have
$$
F(xy)=F(x)F(y), \qquad xV(y)=V(F(x)y)
$$
for any $n$, $i$, $j$, $x \in M^i_n$, $y \in M^j_n$.
\end{enumerate}
\end{defn}

We note that Definition~\ref{FV.pro}~\thetag{ii} immediately yields
$pFd = dF$, $pdV=Vd$. We also note that every $FV$-de Rham
procomplex such that the restriction maps $R$ are surjective and
$$
M_\idot^\hdot = \lim_{\overset{n}{\gets}}M^\hdot_n
$$
has no $p$-torsion is automatically a $V$-procomplex -- for example,
to check Definition~\ref{V.pro}~\thetag{iii}, note that
$$
\begin{aligned}
F(V(y)d\ux) &= FV(y)F(d\ux)=ydF(\ux)=py\ux^{p-1}d\ux\\
&=F(V(y\ux^{p-1}d\ux))=F(V(y\ux^{p-1})V(d\ux)),
\end{aligned}
$$
so that both sides become equal after applying $F$, hence also after
multiplying by $p=FV$.

All $V$-de Rham procomplexes on $X$ form a category in an obvious
way. A major result of \cite{ill} that we will need can be
summarized as follows.

\begin{theorem}[\cite{ill}]\label{ill.main}
\begin{enumerate}
\item The category of $V$-de Rham procomplexes on $X$ admits an
  initial object $W_\idot\Omega^\hdot_A$.
\item This $V$-de Rham procomplex $W_\idot\Omega_A$ has no
  $p$-torsion and uniquely extends to an $FV$-de Rham procomplex.
\end{enumerate}
\end{theorem}

\proof{} \thetag{i} is \cite[Th\'eor\`eme I.1.3]{ill}.  In
\thetag{ii}, lack of $p$-torsion is \cite[Corollaire I.3.5]{ill},
and existence is \cite[Th\'eor\`eme I.2.17]{ill} and
\cite[Proposition I.2.18]{ill}. Uniqueness immediately follows from
the lack of $p$-torsion, since $VF=p$.
\endproof

\begin{defn}\label{dr.witt.def}
The $FV$-de Rham procomplex $W_\idot\Omega^\hdot_A$ of
Theorem~\ref{ill.main}~\thetag{ii} is the {\em de Rham-Witt complex}
of the scheme $X=\Spec A$.
\end{defn}

\begin{remark}
In fact, \cite{ill} works for any smooth finite-type scheme $X/k$,
not necessarily an affine one; we will not need this.
\end{remark}

\subsection{Normalization.}

Next, we need a convenient way to check whether a given $FV$-de Rham
procomplex is the de Rham-Witt complex or not.

\medskip

Recall that the de Rham complex $\Omega^\hdot_A$ has the following
universal property: for any commutative DG algebra $M^\hdot$ over
$k$, any algebra map $A \to M^0$ extends uniquely to a DG algebra
map $\Omega^\hdot_A \to M^\hdot$. In particular, for any $FV$-de
Rham procomplex $M^\hdot_\idot$, the isomorphism $M_1^0 \cong W_1(A)
= A$ induces a map
\begin{equation}\label{omega.m.1}
\Omega^\hdot_A \to M^\hdot_1.
\end{equation}
As in Subsection~\ref{car.subs}, let
$B\Omega^\hdot_A,Z\Omega^\hdot_A \subset \Omega^\hdot_A$ be the
subcomplexes of exact resp.\ closed forms. Recall that we have the
Cartier isomorphism $C$ of \eqref{car.c}; composing $C$ and its
inverse with the natural embedding and natural projection, we obtain
canonical maps
\begin{equation}\label{car.inv}
\wt{C}:Z\Omega^\hdot_A \to \Omega^\hdot_A, \qquad \wt{C}':\Omega_A \to
\Omega^\hdot_A/B\Omega^\hdot_A,
\end{equation}
where from now on, we will ignore $k$-vector space structures and
Frobenius twists. By induction, denote $B_1\Omega^\hdot_A =
B\Omega^\hdot_A$, $Z_1\Omega^\hdot_A = Z\Omega^\hdot_A$ and
\begin{equation}\label{B.n.Z.n}
B_n\Omega^\hdot_A = \wt{C}^{-1}(B_{n-1}\Omega^\hdot_A),\ 
Z_n\Omega^\hdot_A = \wt{C}^{-1}(Z_{n-1}\Omega^\hdot_A) \subset
Z\Omega_A^\hdot,
\end{equation}
for any $n \geq 2$. Assume given an $FV$-de Rham procomplex
$M^\hdot_\idot$ on $X$, and define $\ol{M}^\hdot_\idot$ by
\begin{equation}\label{ol.M}
\ol{M}^\hdot_n = \Ker R \subset M^\hdot_n
\end{equation}
for any $n \geq 1$. Note that $\ol{M}^\hdot_1=M^\hdot_1$, and since
$V$ commutes with $R$, the natural map $V^n$ induces a map
\begin{equation}\label{ol.V}
\ol{V^n}:M_1^\hdot \to \ol{M}_n^\hdot
\end{equation}
for any $n \geq 1$.

\begin{lemma}\label{R.surj}
Assume given an $FV$-de Rham procomplex $M^\hdot_\idot$ such that
the natural map \eqref{omega.m.1} is an isomorphisms, and for every
$i$ and $n' < n$ with some fixed $n$, we have
\begin{equation}\label{V.spans}
\ol{M}^i_n = \Im \ol{V^n} + d \Im \ol{V^n}.
\end{equation}
Then the restriction maps $R:M^\hdot_{n'+1} \to M^\hdot_{n'}$ are
surjective for any $n' < n$.
\end{lemma}

\proof{} For every $n$, any element $x \in M^0_1$ can be lifted to
an element of $M^0_n = W_n(M^0_1)$ (for example, to the
Teichm\"uller representative $\ux$). Thus $R:M^0_n \to M^0_1$ is
surjective. Since $\Omega^i_X$ is spanned by forms $x_0dx_1 \wedge
\dots \wedge dx_i$, $R:M^i_n \to M^i_1$ is surjective for every
$i$. To deduce that $R:M^i_n \to M^i_{n'}$ is surjective for every
$n'<n$, use induction on $n'$ and \eqref{V.spans}.
\endproof

\begin{defn}\label{FV.norm}
An $FV$-de Rham procomplex $M^\hdot_\idot$ is {\em normalized} if
the following holds.
\begin{enumerate}
\item The natural map \eqref{omega.m.1} is an isomorphism.
\item We have a commutative diagram
$$
\begin{CD}
M^\hdot_2 @>{F}>> M_1^\hdot \cong \Omega^\hdot_A\\
@V{R}VV @VV{\pi}V\\
M^\hdot_1 \cong \Omega^\hdot_A @>{\wt{C}'}>>
\Omega^\hdot_A/B\Omega^\hdot_A,
\end{CD}
$$
where $\pi$ is the natural projection map, and $\wt{C}'$ is the map
\eqref{car.inv}.
\item For any $n \geq 1$, the map
$$
\eta = \ol{V^n} \oplus (d \circ \ol{V^n}): M^\hdot_1 \oplus
M^{\hdot-1}_1 \to \ol{M}^\hdot_n
$$
is surjective, and its kernel $\Ker \eta \subset \Omega^\hdot_A
\oplus \Omega^{\hdot-1}_A = M^\hdot_1 \oplus M^{\hdot-1}_1$ is an
extension of $Z_n\Omega^{\hdot-1}_A \subset \Omega^{\hdot-1}_A$ by
$B_n\Omega^\hdot_A \subset \Omega^\hdot_A$ --- that is, we have a
commutative diagram
$$
\begin{CD}
0 @>>> B_n\Omega^\hdot_A @>>> \Ker \eta @>>> Z_n\Omega^{\hdot-1}_A
@>>> 0\\
@. @VVV @VVV @VVV\\
0 @>>> \Omega^\hdot_A @>>> \Omega^\hdot_A \oplus \Omega^{\hdot-1}_A
@>>> \Omega^{\hdot-1}_A @>>> 0
\end{CD}
$$
with exact rows.
\end{enumerate}
\end{defn}

\begin{prop}\label{crit.prop}
An $FV$-de Rham procomplex is normalized if and only if it is
isomorphic to the de Rham-Witt complex $W_\idot\Omega^\hdot_A$.
\end{prop}

The difficult part of this statement is the ``if'' part, but
fortunately, this has been proved in \cite{ill}. To prove the ``only
if'' part, we need the following lemma.

\begin{lemma}\label{no.p}
If an $FV$-de Rham procomplex $M^\hdot_\idot$ is normalized, then
$$
M^i_\idot = \lim_{\overset{n}{\gets}}M^i_n
$$
has no $p$-torsion for any $i$.
\end{lemma}

\proof{} Since we have $p=FV$, all the groups $\ol{M}^i_n$ are
annihilated by $p$, so that multiplication by $p$ induces a map
$$
\ol{p}_n:\ol{M}^i_n \to \ol{M}^i_{n+1}.
$$
It suffices to prove that this map is injective for every $n \geq
1$. Moreover, we have $\ol{p}_n \circ \ol{V^n} = \ol{V^{n+1}} \circ
F$, so that $\ol{p}_n$ sends $\Im \ol{V^n}$ into $\Im \ol{V^{n+1}}$,
and it suffices to prove that the induced maps
$$
\ol{p}_n:\Im \ol{V^n} \to \Im \ol{V^{n+1}}, \qquad \ol{p}_n:\Coker
\ol{V^n} \to \Coker \ol{V^{n+1}}
$$
are injective. Using Definition~\ref{FV.norm}~\thetag{iii} and
\thetag{iv}, we rewrite these maps as maps 
$$
\Omega^i_A/B_n\Omega^i_A \to \Omega^i_A/B_{n+1}\Omega^i_A, \qquad
\Omega^{i-1}_A/Z_{n-1}\Omega^{i-1}_A \to
\Omega^{i-1}_A/Z_n\Omega^{i-1}_A,
$$
and since $pV^n=V^{n+1}F$ and $pdV^n=V^{n+1}F$, by
Definition~\ref{FV.norm}~\thetag{ii}, both maps are induced by the
inverse Cartier map $\wt{C}'$ of \eqref{car.inv}. Hence both are
injective.
\endproof

\proof[Proof of Proposition~\ref{crit.prop}.] For the ``if'' part,
Definition~\ref{FV.norm}~\thetag{ii} is \cite[Proposition
  I.3.3]{ill}, and the rest is \cite[Corollaire
  I.3.9]{ill}. Conversely, assume given a normalized $FV$-de Rham
procomplex $M^\hdot_\idot$ on $X$. Then by
Definition~\ref{FV.norm}~\thetag{i}, \thetag{iii}, $M^\hdot_\idot$
satisfies the assumptions of Lemma~\ref{R.surj} for any $n \geq 1$,
and by Lemma~\ref{R.surj} and Lemma~\ref{no.p}, forgetting $F$ turns
it into a $V$-de Rham procomplex in the sense of
Definition~\ref{V.pro}. Then by Theorem~\ref{ill.main}, we have
a natural map
$$
\phi:W_\idot\Omega^\hdot_A \to M^\hdot_\idot
$$
from the de Rham-Witt complex $W_\idot\Omega^\hdot_A$. This map is
compatible with the restriction maps, so that it suffices to prove
that for any $n \geq 1$, $\phi$ identifies the kernel of the map
$$
R:W_{n+1}\Omega^\hdot_A \to W_n\Omega^\hdot_A
$$
with $\ol{M}^\hdot_n$. Since the de Rham-Witt complex is normalized,
this immediately follows from Definition~\ref{FV.norm}~\thetag{iii}
and \thetag{iv}.
\endproof

\subsection{Comparison.}

Now keep the assumptions and notation of the last Subsection, and
consider Hochschild-Witt homology groups $W_\idot HH_\idot(A)$ of
the algebra $A$. Since $A$ is commutative, $W_nHH_\idot(A)$ is a
graded-commutative algebra with respect to the product
\eqref{WHH.prod}, and the Connes-Tsygan differential $B$ is a
derivation for this product. Moreover, we have natural maps $V$ and
$F$ of \eqref{FV.WHH}.

\begin{lemma}\label{w.fv.le}
The groups $M^i_n = W_nHH_i(A)$ with the product \eqref{WHH.prod},
the differential $d=B$, and the maps $V$, $F$ of \eqref{FV.WHH} form
an $FV$-de Rham procomplex on $X = \Spec A$ in the sense of
Definition~\ref{FV.pro}.
\end{lemma}

\proof{} Definition~\ref{FV.pro}~\thetag{ii} is
Lemma~\ref{fdv.whh.le}, and \thetag{iii} immediately follows from
Lemma~\ref{FV.prod.le}. It remains to check \thetag{i}. We have the
functorial isomorphisms $W_nHH_0(A) \cong W_n^H(A)$ provided by
Theorem~\ref{hess.thm}, and they are compatible with the maps $R$
and $V$. Combining these with the identifications $W_n^H(A) \cong
W_n^{cl}(A)$, we obtain isomorphisms $\iota:W_nHH_0(A) \cong
W_n^{cl}(A)$, $n \geq 1$ that are also compatible with $R$ and
$V$. In the limit, we get an identification $\iota:WHH_0(A) \cong
W(A)$. Since this group has no $p$-torsion, $FV=VF=p$ implies that
$\iota$ is also compatible with $F$, and then the same is true for
its restrictions $\iota:W_nHH_0(A) \cong W_n^{cl}(A)$. It remains to
show that $\iota$ is multiplicative. But every element $x \in
W_nHH_0(A)$ is of the form $V(y) + T(a)$, $y \in W_{n-1}HH_0(A)$, $a
\in A$, where $T$ is the Teichm\"uller map \eqref{T.A}. Therefore by
induction and \thetag{iii}, it suffices to prove that $\iota$ is
multiplicative on the image of $T$. But $T$ is multiplicative by
\eqref{T.mult}, and $\iota(T(a))$, $a \in A$ is the Teichm\"uller
representative $\underline{a}$ of $a$.
\endproof

\begin{lemma}\label{FV.ii}
The $FV$-de Rham procomplex $W_\idot HH_\idot(A)$ of
Lemma~\ref{w.fv.le} satisfies the condition \thetag{ii} of
Definition~\ref{FV.norm}.
\end{lemma}

\proof{} By definition, in terms of the identifications of
Lemma~\ref{la.ho} and \eqref{K.l.cor}, the Frobenius map
$F:W_2HH_\idot(A) \to HH_\idot(A)$ is induced by the composition map
$$
\begin{CD}
\K_\idot(W_2A^\hush) @>{\K_\idot(\overline{F})}>>
\K_\idot(\pi_{p*}i_p^*A^\hush) @>{\phi_p}>> \K^p_\idot(A^\hush),
\end{CD}
$$
where $\phi_p$ is the map \eqref{nu.phi.eq}. By Lemma~\ref{nu.phi},
this map coincides with the map $\e_p$ on homology of degree $0$,
and by \cite[Lemma~4.1~\thetag{iii}]{cart}, we have $\e_p = C \circ
R$. Therefore $\phi_p$ takes values in the subcomplex
$Z\K^p_\idot(A^\hush) \subset K^p_\idot(A^\hush)$. Moreover, we have
a tautological commutative diagram
$$
\begin{CD}
Z\K^p(A^\hush) @>>> \K^p(A^\hush)\\
@V{\xi}VV @VVV\\
Z\K^p(A^\hush)/B\K^p(A^\hush) @>>> \K^p(A^\hush)/B\K^p(A^\hush),
\end{CD}
$$
where $\xi$ is the natural projection, so that to prove the claim,
it suffices to prove that we have a commutative diagram
\begin{equation}\label{w2.dia}
\begin{CD}
\K_\idot(W_2A^\hush) @>{\K_\idot(\ol{F})}>> \K_\idot(\pi_{p*}i_p^*A^\hush)\\
@V{\K_\idot(R)}VV @VV{\xi \circ \phi_p}V\\
\K_\idot(A^\hush) @>{\sim}>> Z\K^p_\idot(A^\hush)/B\K^p_\idot(A^\hush)
\end{CD}
\end{equation}
in $\D(\Lambda,W_2(k))$, where the map in the bottom is the
isomorphism \eqref{car.nc}. As in the proof of
Theorem~\ref{car.thm}~\thetag{i} given in
Subsection~\ref{basic.subs}, let $E = i_p^*A^\hush$, and consider
the complexes $E'_\idot$, $\wt{E}_\idot$. Then bottom map in
\eqref{w2.dia} is the map \eqref{bb.eq} associated to the map
$b:E'_\idot \to \wt{E}_\idot$, and the rightmost map is the map
\eqref{bb.eq} associated to the tautological embedding $\pi_{p*}E
\cong \wt{E}_0 \subset \wt{E}_\idot$. Thus to prove that the diagram
is commutative, it suffices to observe that we have a commutative
diagram
$$
\begin{CD}
W_2A^\hush \cong E'_0 @>{\overline{F}}>> \pi_{p*}E \cong \wt{E}_0\\
@VVV @VVV\\
E'_\idot @>{b}>> \wt{E}_\idot,
\end{CD}
$$
where the vertical arrows are the natural embedding maps, and the map
$$
\begin{CD}
W_2A^\hush @>>> E'_\idot @>>> A^\hush
\end{CD}
$$
obtained by composing the embedding $E'_0 \subset E'_\idot$ with the
quasiisomorphism \eqref{comp.bb} is the restriction map $R$.
\endproof

To proceed further, we need to adapt \eqref{B.Z.K} to describe the
subgroups $B_n\Omega^\hdot_A$, $Z_n\Omega^\hdot_A$ of
\eqref{B.n.Z.n}. To this effect, consider the iterates
\begin{equation}\label{psi.n}
C^n:A^\hush \to \pi_{p^n!}i_{p^n}^*A^\hush, \qquad
R^n:\pi_{p^n*}i_{p^n}^*A^\hush \to A^\hush.
\end{equation}
of the maps \eqref{C.eq}, and simplify notation by setting
$\K^{(n)}_\idot(A^\hush) = \K^{p^n}_\idot(A^\hush)$ for any $n \geq
1$. Then \eqref{K.l.4} induces natural augmentation maps
$$
\begin{CD}
\pi_{p^n*}i_{p^n}^*A^\hush[1] @>{\wt{\kappa}_1}>>
\K^{(n)}(A^\hush) @>{\wt{\kappa}_0}>>
\pi_{p^n!}i_{p^n}^*A^\hush
\end{CD}
$$
and we can define subcomplexes $B\K^{(n)}_\idot(A^\hush),
Z\K^{(n)}_\idot(A^\hush) \subset \K^{(n)}_\idot(A^\hush)$ by
\begin{equation}\label{B.Z.K.n}
B\K^{(n)}_\idot(A^\hush) = \wt{\kappa}_1(\Ker R^n), \qquad
Z\K^{(n)}_\idot(A^\hush) = \wt{\kappa}_0^{-1}(\Im C^n).
\end{equation}
Since $\wt{\kappa}_0 \circ \wt{\kappa}_1 = 0$, we have
$B\K^{(n)}_\idot(A^\hush) \subset Z\K^{(n)}(A^\hush)$. Moreover,
\eqref{K.l.cor} provides a canonical identification
\begin{equation}\label{K.(n)}
HC_\idot(\K^{(n)}_\idot(A_\hush))
\cong HH_\idot(A),
\end{equation}
and if we let
$$
B_nHH_\idot(A) = HC_\idot(BK^{(n)}_\idot(A^\hush)), \quad
Z_nHH_\idot(A) = HC_\idot(ZK^{(n)}_\idot(A^\hush)),
$$
then we have natural maps
\begin{equation}\label{b.z.n}
B_nHH_\idot(A) \to Z_nHH_\idot(A) \to HH_\idot(A).
\end{equation}

\begin{lemma}\label{b.z.n.le}
For any integer $n \geq 1$, the natural maps \eqref{b.z.n} are
injective, and if we identify $HH_\idot(A) \cong \Omega^\hdot_A$ by
Theorem~\ref{hkr}, then $B_nHH_\idot(A) \cong B_n\Omega^\hdot_A
\subset \Omega^\hdot_A$ and $Z_nHH_\idot(A) \cong Z_n\Omega^\hdot_A
\subset \Omega^\hdot_A$.
\end{lemma}

\proof{} Induction on $n$. For $n=1$, \eqref{B.Z.K.n} coincides with
\eqref{B.Z.K}, so that we are done by
Theorem~\ref{car.thm}~\thetag{ii}. For any $n \geq 2$,
Lemma~\ref{free.le} shows that $\pi_{p^{n-1}!}i_{p^{n-1}}^*\Phi
A^\hush \cong \pi_{p^{n-1}*}i_{p^{n-1}}^*\Phi A^\hush$, and
moreover, if we denote this object by $\Phi^{(n)}A^\hush$, then we
have natural short exact sequences
$$
\begin{CD}
0 @>>> \pi_{p^{n-1}!}i_{p^{n-1}}^*A^\hush @>{C}>>
\pi_{p^n!}i_{p^n}^*A^\hush @>>> \Phi^{(n)}A^\hush @>>> 0,\\
0 @>>> \Phi^{(n)}A^\hush @>>> \pi_{p^n*}i_{p^n}^*A^\hush
@>{R}>> \pi_{p^{n-1}*}i_{p^{n-1}}^*A^\hush @>>> 0.
\end{CD}
$$
Furthermore, we have $\K^{(n)}_\idot(A^\hush) \cong
\pi_{p^{n-1}!}i_{p^{n-1}}^*\K^{(1)}_\idot(A^\hush)$ by base change,
and if we denote
$$
\ol{B}\K^{(n)}_\idot(A^\hush) = \wt{\kappa}_1(\Ker
R), \qquad \ol{Z}\K^{(n)}_\idot(A^\hush) = \wt{\kappa}_0^{-1}(\Im
C),
$$
then we have
$$
\ol{B}\K^{(n)}_\idot(A^\hush) \cong \Phi^{(n)}A^\hush[-1] \cong
L^\hdot\pi_{p^{n-1}!}i_{p^{n-1}}^*B\K^{(1)}_\idot(A^\hush)
$$
and
$$
\ol{Z}\K^{(n)}_\idot(A^\hush) \cong
L^\hdot\pi_{p^{n-1}!}i_{p^{n-1}}^*Z\K^{(1)}_\idot(A^\hush).
$$
Therefore when we pass to homology, the maps
$$
HC_\idot(\ol{B}\K^{(n)}_\idot(A^\hush)) \to
HC_\idot(\ol{Z}\K^{(n)}_\idot(A^\hush)) \to
HC_\idot(\K^{(n)}_\idot(A^\hush))
$$
are injective, and under \eqref{K.(n)}, these are exactly the maps
\eqref{b.z.n} for $n=1$. Moreover, possibly after applying the
tautological functor $\nu$, the isomorphism \eqref{car.nc} induces a
natural quasiisomorphism
\begin{equation}\label{car.n}
\overline{Z}\K^{(n)}_\idot(A^\hush)/\overline{B}\K^{(n)}(A^\hush)
\cong \K^{(n-1)}_\idot(A^\hush)
\end{equation}
that gives the inverse Cartier map under the identification
\eqref{K.(n)}. It remains to notice that by definition, we have
$$
\overline{B}\K^{(n)}_\idot(A^\hush) \subset B\K^{(n)}_\idot(A^\hush)
\subset Z\K^{(n)}_\idot(A^\hush) \subset
\overline{Z}\K^{(n)}_\idot(A^\hush),
$$
and the quasiisomorphism \eqref{car.n} provides exact triangles
$$
\begin{CD}
\overline{B}\K^{(n)}_\idot(A^\hush) @>>> B\K^{(n)}_\idot(A^\hush)
@>>> B\K^{(n-1)}_\idot(A^\hush) @>>>\\
\overline{B}\K^{(n)}_\idot(A^\hush) @>>> Z\K^{(n)}_\idot(A^\hush)
@>>> Z\K^{(n-1)}_\idot(A^\hush) @>>>
\end{CD}
$$
Comparing this with \eqref{B.n.Z.n} and applying induction, we get
the claim.
\endproof

\begin{remark}
Note that \eqref{B.Z.K.n} makes sense for any associative unital
$k$-algebra $A$, so that the maps \eqref{b.z.n} also exist in full
generality. Moreover, for any $m \leq n$, one can consider the
subcomplexes
$$
B_m\K^{(n)}_\idot(A^\hush) = \wt{\kappa}_1(\Ker R^m), \qquad
Z_m\K^{(n)}_\idot(A^\hush) = \wt{\kappa}_0^{-1}(\Im C^m),
$$
and then the same argument as in the proof of Lemma~\ref{b.z.n.le}
shows that we have natural identifications
$$
B_mHH_\idot(A) \cong HC_\idot(B_m\K^{(n)}_\idot(A^\hush)),
\quad
Z_mHH_\idot(A) \cong HC_\idot(Z_m\K^{(n)}_\idot(A^\hush)).
$$
Thus for any associative unital $k$-algebra $A$, we have an infinite
sequence of natural maps
\begin{equation}\label{b.n.z.n.a}
\begin{aligned}
B_1HH_\idot(A) &\to \dots \to B_nHH_\idot(A) \to \dots\\
\dots &\to Z_nHH_\idot(A) \to \dots \to Z_1HH_\idot(A) \to HH_\idot(A),
\end{aligned}
\end{equation}
and for any $n$, the first $n$ terms in the sequence are represented
explicitly by subcomplexes in $\K^{(n)}_\idot(A^\hush)$. We also
have natural long exact sequences
$$
\begin{CD}
B_nHH_\idot(A) @>>> Z_nHH_\idot(A) @>>> HH_\idot(A) @>>>
\end{CD}
$$
obtained by iterating the quasiisomorphism \eqref{car.nc}. We do not
know under what assumptions the maps \eqref{b.n.z.n.a} are
injective.
\end{remark}

\begin{theorem}\label{ill.thm}
Assume given a $k$-algebra $A$ satisfying the assumptions of
Theorem~\ref{hkr}, and let $W_\idot\Omega^\hdot_A$ be the de
Rham-Witt complex of the scheme $X = \Spec A$, as in
Definition~\ref{dr.witt.def}. Then the Hochschild-Kostant-Rosenberg
isomorphism \eqref{hkr.iso} extends to a series of functorial
multiplicative isomorphisms
$$
W_nHH_i(A) \cong W_n\Omega^i_A
$$
that commute with the maps $V$, $F$, and send the Connes-Tsygan
differential $B$ to the de Rham-Witt differential $d$.
\end{theorem}

\proof{} By Proposition~\ref{crit.prop}, it suffices to prove that
the $FV$-de Rham procomplex $M^\hdot_\idot = W_\idot HH_\idot(A)$ of
Lemma~\ref{w.fv.le} is normalized in the sense of
Definition~\ref{FV.norm}. Definition~\ref{FV.norm}~\thetag{i} is
Theorem~\ref{B.d.thm}, and Definition~\ref{FV.norm}~\thetag{ii} is
Lemma~\ref{FV.ii}, so what we have to check is
Definition~\ref{FV.norm}~\thetag{iii}. Use induction on $n$. The
base case $n=0$ is trivial. Assume the statement proved for all $n'
< n$. The first of the exact sequences of Proposition~\ref{whh.prop}
induces a long exact sequence
$$
\begin{CD}
HH_\idot(\pi_{p^n!}i_{p^n}^*(A^\hush)) @>>>
W_nHH_\idot(A) @>{R}>> W_{n-1}HH_\idot(A) @>>>
\end{CD}
$$
of homology groups. Note that Definition~\ref{FV.norm}~\thetag{iii}
implies \eqref{V.spans}, so that in particular, $M_{n'}^\hdot$
satisfy the assumptions of Lemma~\ref{R.surj}. Then $R$ is
surjective, so that the connecting differential in the long exact
sequence vanishes, and we have
$$
\ol{M}^\hdot_n \cong HH_\idot(\pi_{p^n!}i_{p^n}^*A^\hush) \cong
HC_\idot(\K_\idot(\pi_{p^n!}i_{p^n}^*A^\hush)).
$$
By definition, the map $\ol{V^n}$ is then induced by the canonical map
\begin{equation}\label{nu.pn}
\nu_{p^n}:\K^{(n)}(A^\hush) \to
\K_\idot(\pi_{p^n!}i_{p^n}^*A^\hush)
\end{equation}
of \eqref{nu.phi.eq}, where we have used the identification
\eqref{K.(n)}. By Lemma~\ref{nu.phi}, this map is an isomorphism in
homological degree $0$, and equal to the canonical map $\e_{p^n}$ in
homological degree $1$. By \cite[Lemma~4.1~\thetag{iii}]{cart}, we
have $\e_{p^n} = C^n \circ R^n$, where $C^n$ and $R^n$ are the maps
\eqref{psi.n}. In terms of the complexes \eqref{B.Z.K.n}, this means
that we have $\Ker \nu_{p^n} \cong B\K^{(n)}_\idot(A^\hush) \subset
\K^{(n)}_\idot(A^\hush)$, and the cokernel $\Coker \nu_{p^n}$ of the
map \eqref{nu.pn} fits into a short exact sequence
$$
\begin{CD}
0 @>>> Z\K^{(n)}_\idot(A^\hush)[1] @>>> \K^{(n)}_\idot(A^\hush)[1]
@>{B\circ\nu_{p^n}}>> \Coker \nu_{p^n} @>>> 0,
\end{CD}
$$
where $B:\K_\idot(\pi_{p^n!}i_{p^n}^*A^\hush)[1] \to
\K_\idot(\pi_{p^n!}i_{p^n}^*A^\hush)$ is the map \eqref{B.gen}. We
conclude that the map
\begin{equation}\label{nu}
\eta=\nu_{p^n} \oplus (B \circ \nu_{p^n}):\K^{(n)}_\idot(A^\hush)
\oplus \K^{(n)}_\idot(A^\hush)[1] \to
\K_\idot(\pi_{p^n!}i_{p^n}^*A^\hush)
\end{equation}
is surjective, and its kernel is an extension of
$Z\K^{(n)}_\idot(A^\hush)[1]$ by $B\K^{(n)}_\idot(A^\hush)$ in the
same sense as in Definition~\ref{FV.norm}~\thetag{iii}. Therefore we
have a commutative diagram
$$
\begin{CD}
B_nHH_\idot(A) @>>> HC_\idot(\Ker \eta) @>>> Z_nHH_{\idot-1}(A)\\
@VVV @VVV @VVV\\
HH_\idot(A) @>>> HH_\idot(A) \oplus HH_{\idot-1}(A) @>>>
HH_{\idot-1}(A)
\end{CD}
$$
whose rows are exact in the middle term. By
Theorem~\ref{car.thm}~\thetag{iii}, the leftmost and the rightmost
vertical maps in this diagram are injective. Since the top row comes
from a long exact sequence, and the bottom row is exact on the left
and the right, the top row is then also exact on the left and on the
right. Moreover, the middle vertical map is also injective, so that
the connecting differential
$$
HC_\idot(\K_\idot(\pi_{p^n!}i_{p^n}^*A^\hush)) \cong HC_\idot(\Coker
\eta) \to HC_{\idot-1}(\Ker \eta)
$$
in the long exact sequence induced by the surjective map \eqref{nu}
vanishes. Therefore the map
$$
\begin{aligned}
HH_\idot(A) \oplus HH_{\idot-1}(A) &\cong
HC_\idot(\K^{(n)}_\idot(A^\hush) \oplus \K^{(n)}_\idot(A^\hush)[1])
\to\\
&\to HC_\idot(\K_\idot(\pi_{p^n!}i_{p^n}^*A^\hush)) \cong \ol{M}^\hdot_n
\end{aligned}
$$
induced by the map \eqref{nu} is surjective, and its kernel
$HC_\idot(\Ker \eta)$ is an extension of $B_nHH_\idot(A)$ by
$Z_nHH_{\idot-1}(A)$. This is Definition~\ref{FV.norm}~\thetag{iii}.
\endproof

{\small\noindent
Affiliations (in the precise form required for legal reasons):
\begin{enumerate}
\renewcommand{\labelenumi}{\arabic{enumi}.}
\item Steklov Mathematics Institute, Algebraic Geometry Section
  (main affiliation).
\item Laboratory of Algebraic Geometry, National Research University
Higher\\ School of Economics.
\item Center for Geometry and Physics, Institute for Basic
  Science (IBS), Pohang, Korea.
\end{enumerate}}

\noindent
{\em E-mail address\/}: {\tt kaledin@mi.ras.ru}

\end{document}